\documentclass[10pt]{article}
 \usepackage[dvipdfm]{graphicx, color}
 \usepackage{amssymb, amsmath,amsthm}
\newtheorem{thm}{{\bf Theorem}}[section]
 \newtheorem{lem}[thm]{{\bf Lemma}}
  \newtheorem{cor}[thm]{{\bf Corollary}}
  \newtheorem{prop}[thm]{{\bf Proposition}}

  \newtheorem{rem}[thm]{Remark}

  \newtheorem{ques}[thm]{Question}

   \numberwithin{equation}{section}

\begin{document}
\title{Pseudo-Anosovs on closed surfaces having small entropy and the Whitehead sister link exterior}
 \author{Eiko Kin\footnote{The author is partially supported by Grant-in-Aid for Young Scientists (B) (No. 20740031), 
MEXT, Japan. 
\hspace{3cm}
Last compile 08/15/2011. 
}\hspace{1mm}
and Mitsuhiko Takasawa
}
\date{\empty}
\maketitle

\begin{abstract}
We denote by $\delta_g$ (resp. $\delta_g^+$),  the minimal dilatation for pseudo-Anosovs 
(resp. pseudo-Anosovs with orientable invariant foliations)  on a closed surface  of genus $g$. 
This paper concerns the  pseudo-Anosovs  which occur as  monodromies of fibrations on manifolds 
obtained from the Whitehead sister link exterior $W$ by Dehn filling  two cusps, 
where the fillings are on the boundary slopes of fibers of $W$. 
We give upper bounds of $\delta_g$ for $g \equiv 0,1,5,6,7,9 \pmod{10}$, 
$\delta_g^+$ for $g \equiv 1,5,7,9 \pmod{10}$. 
Our bounds improve the previous one given by Hironaka. 
We note that the monodromies of fibrations on $W$ were also studied by Aaber and Dunfield independently. 
\end{abstract}

\section{Introduction}

Let $ \mathrm{Mod}(\varSigma)$ be the mapping class group on an orientable surface $\varSigma$. 
An element $\phi \in \mathrm{Mod}(\varSigma)$ which contains a pseudo-Anosov homeomorphism $\Phi: \varSigma \rightarrow \varSigma$ 
as a representative  is called a pseudo-Anosov mapping class.  
There are two numerical invariants for  pseudo-Anosov mapping classes. 
One  is the dilatation $\lambda(\phi)>1$ (or the entropy $\mathrm{ent}(\phi) = \log \lambda(\phi)$) which is defined to be the dilatation $\lambda(\Phi)$ of $\Phi$, and 
the other  is the hyperbolic volume $\mathrm{vol}(\phi)= \mathrm{vol}({\Bbb T}(\phi))$  of the mapping torus  ${\Bbb T}(\Phi)$.  
It is natural to ask whether there is a relation between $\mathrm{ent}(\phi)$ and $\mathrm{vol}(\phi)$. 
Computer experiments in \cite{KKT} tell us that 
if we fix a surface $\varSigma$, then 
pseudo-Anosovs with small dilatation have small volume. 
This is true in a sense. 
In fact  it is proved in \cite{FLM} that pseudo-Anosovs on any surfaces with small dilatation have the bounded volume, see Theorem~\ref{thm_finite-FLM}.

We denote by $\delta_g$, the minimal dilatation for  pseudo-Anosov elements $\phi \in  \mathrm{Mod}(\varSigma_g)$ 
on  a closed surface $\varSigma_g$ of genus $g$. 
A natural question is: what is  the value $\delta_g$? 
To discuss the minimal dilatations, we introduce the polynomial  
$$ f_{(k,\ell)}(t)= t^{2k} - t^{k+\ell} - t^k - t^{k-\ell}+1 \hspace{2mm}\mbox{for}\  k >0,\ -k < \ell < k.$$
This polynomial has the largest real root  $\lambda_{(k,\ell)}$ which is greater than $1$ (Theorem~\ref{thm_Fried-Oertel-Poly} and Lemma~\ref{lem_dilatation_kl}). 
For any fixed $\ell>0$, it follows that  
$k  \log \lambda_{(k,\ell)} $  converges to $\log  (\tfrac{3 + \sqrt{5}}{2})$ if $k $ goes to $\infty$ (Lemma~\ref{lem_asymp_roots}). 
It is easy to show that 
$\delta_1= \lambda_{(1,0)}= \tfrac{3+ \sqrt{5}}{2}$. 
It was proved by Cho-Ham that 
$\delta_2 = \lambda_{(2,1)} \approx 1.72208$ \cite{CH}. 
It is open to determine the values $\delta_g$ for $g \ge 3$. 
Questions on properties of $\delta_g$ were posed by McMullen and Farb: 

\begin{ques}[\cite{McMullen} for (1), \cite{Farb} for (2)]
\label{ques_Mc-Farb}
\begin{description}
\item[(1)] 
Does $\displaystyle \lim_{g \to \infty}g  \log \delta_g$ exist? What is its value?
\item[(2)] 
Is the sequence $\{\delta_g\}_{g \ge 2}$ (strictly) monotone decreasing? 
\end{description}
\end{ques}

\noindent
Related questions  are ones  for orientable pseudo-Anosovs.  
A pseudo-Anosov mapping class $\phi$ is said to  be {\it orientable} 
if the invariant (un)stable foliation of a pseudo-Anosov homeomorphism $\Phi \in  \phi$ is orientable. 
We denote by $\delta_g^+$, the minimal dilatation for orientable pseudo-Anosov elements of $ \mathrm{Mod}(\varSigma_g)$. 
%Obviously $\delta_{g,p} \le \delta_{g,p}^+$. 
The minima $\delta_g^+$ were determined for $g=2$ by Zhirov \cite{Zhirov}, 
for $3 \le g \le 5$ by Lanneau-Thiffeault  \cite{LT}, and  for $g=8$ by Lanneau-Thiffeault and Hironaka \cite{LT,Hironaka}. 
Those values are given by 
$\delta_2^+= \lambda_{(2,1)}$, 
$\delta_3^+ =\lambda_{(3,1)}= \lambda_{(4,3)}  \thickapprox 1.40127 $, 
$\delta_4^+ = \lambda_{(4,1)}  \thickapprox 1.28064$, 
$\delta_5^+ =\lambda_{(6,1)} = \lambda_{(7,4)} \thickapprox 1.17628$ and 
$\delta_8^+=\lambda_{(8,1)} \thickapprox 1.12876$. 

Lanneau-Thiffeault obtained the inequality  $\delta_5^+ \le \delta_6^+$ (\cite{LT})  which implies that 
$\{\delta^+_g\}_{g \ge 2}$ is not strictly monotone decreasing. 
This leads us to ask an alternative question related to Question~\ref{ques_Mc-Farb}(2):  
is the sequence $\{\delta_g^+\}_{g \ge 2}$  monotone decreasing? 
Also, one can ask: 
which $g$ does the inequality $\delta_g < \delta^+_g$ hold? 
It is easy to see that $\delta_1= \delta_1^+$. 
The equality $\delta_g= \delta_g^+$ holds for $g=2$ \cite{CH,Zhirov}. 
We do not know whether $\delta_3=\delta_3^+$ holds or not. 
By work of Lanneau-Thiffeault and Hironaka, it follows that $\delta_g < \delta_g^+$ for $g= 4,6,8$ \cite{LT,Hironaka}.

To discuss Question~\ref{ques_Mc-Farb}(1), 
we recall the previous upper bound of $\delta_g$ given by Hironaka. 
 
\begin{thm}[\cite{Hironaka}]
\label{thm_bound_H}
\begin{description}
\item[(1)] 
$\delta_g \le \lambda_{(g+1,3)}$ if $g \equiv 0,1,3,4 \pmod 6$ and $g \ge 3$. 
\item[(2)] 
$\delta_g \le \lambda_{(g+1,1)}$ if $g \equiv 2,5 \pmod 6$ and $g \ge 5$. 
\end{description}
\end{thm}

\noindent
By using Lemma~\ref{lem_asymp_roots} and Theorem~\ref{thm_bound_H}, the following asymptotic inequality holds. 

\begin{thm}[\cite{Hironaka}]
\label{thm_asymptotic}
$\displaystyle \lim_{g \to \infty} \sup g \log  \delta_g \le \log(\tfrac{3+ \sqrt{5}}{2})$. 
\end{thm}

\noindent
This improves the upper bound $g  \log \delta_g \le  g  \log \delta^+_g \le \log(2 + \sqrt{3})$ for each $g \ge 2$ by Minakawa \cite{Minakawa} and Hironaka-Kin \cite{HK}. 
Since $\log \delta_g$ tends to $0$ as $g $ tends to $\infty$, Theorem~\ref{thm_asymptotic} implies that 
$$ \lim_{g \to \infty} \sup |\chi (\varSigma_g)| \log  \delta_g \le 2 \log(\tfrac{3+ \sqrt{5}}{2}),$$
where $\chi(\varSigma)$ is the Euler characteristic of a surface $\varSigma$.

\begin{figure}
\begin{center}
\includegraphics[width=3.5in]{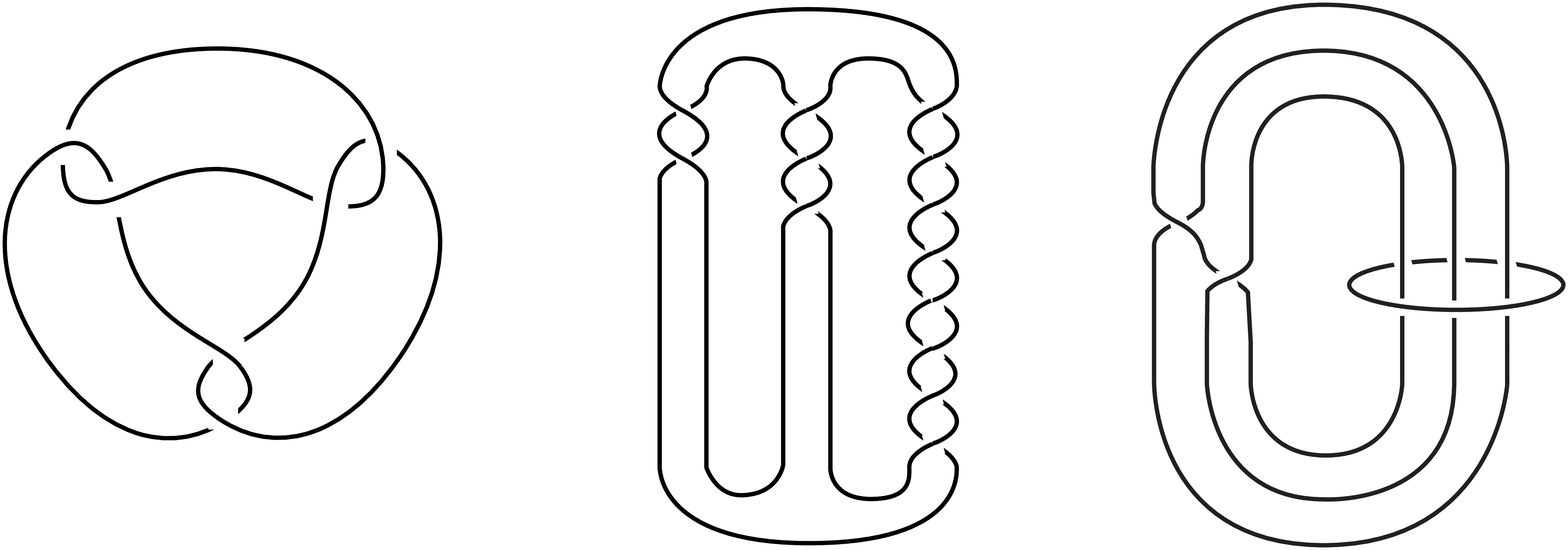}
\caption{(left) $3$ chain link $\mathcal{C}_3$. (center) $(-2,3,8)$-pretzel link or Whitehead sister link. 
(right) link $6_2^2$.}
\label{fig_3chain}
\end{center}
\end{figure}

Let $N$ be the {\it magic manifold} which is the exterior of the $3$ chain link $\mathcal{C}_3$  illustrated in Figure~\ref{fig_3chain}(left).  
This manifold has the smallest known volume among orientable hyperbolic $3$-manifolds having $3$ cusps. 
Many manifolds having at most $2$ cusps with small volume are obtained from $N$ by Dehn fillings, see \cite{MP}.  
In this paper, we study the {\it small dilatation pseudo-Anosov homeomorphisms} which occur as monodromies of fibrations on manifolds 
obtained from $N$ by Dehn filling all three cusps. 
In \cite{FLM}, Farb, Leininger and Margalit introduced small dilatation pseudo-Anosov homeomorphisms which we recall below.

For any number $P >1$,  define the set of pseudo-Anosov homeomorphisms  
$$\Psi_P = \{\mbox{pseudo-Anosov\ }\Phi: \varSigma \rightarrow \varSigma \ |\ \chi(\varSigma)<0,\ |\chi(\varSigma)| \log \lambda(\Phi) \le \log P\}.$$
They call elements $\Phi \in \Psi_P$   small dilatation pseudo-Anosov homeomorphisms. 
Theorem~\ref{thm_asymptotic} says that 
if one takes $P$ sufficiently large, then 
$\Psi_P$ contains a pseudo-Anosov homeomorphism $\Phi_g: \varSigma_g \rightarrow \varSigma_g$ for each $g \ge 2$. 
By a result  by Hironaka-Kin \cite{HK}, $\Psi_P$ also contains a pseudo-Anosov homeomorphism $\Phi_n: D_n \rightarrow D_n$ on an $n$-punctured disk $D_n$ for each $n \ge 3$. 
Let $ \varSigma^{\circ} \subset \varSigma$ be the surface obtained by removing the singularities of the (un)stable foliation for $\Phi$ and 
$\Phi|_{\varSigma^{\circ}}: \varSigma^{\circ} \rightarrow \varSigma^{\circ}$ denotes the restriction. 
Observe that  $\lambda(\Phi)= \lambda(\Phi|_{\varSigma^{\circ}})$. 
The set 
$$\Psi_P^{\circ}= \{\Phi|_{\varSigma^{\circ}}: \varSigma^{\circ} \rightarrow \varSigma^{\circ}\ |\ (\Phi: \varSigma \rightarrow \varSigma) \in \Psi_P\}$$ 
is  infinite. 
Let $\mathcal{T}(\Psi_P^{\circ})$ be the set of homeomorphism classes of mapping tori by elements of $\Psi_P^{\circ}$.

\begin{thm}[\cite{FLM}]
\label{thm_finite-FLM}
The set $\mathcal{T}(\Psi_P^{\circ})$ is  finite. 
Namely, for each $P >1$, there exist finite many complete, non compact hyperbolic $3$-manifolds 
$M_1, M_2, \cdots, M_r$ fibering over $S^1$ so that the following holds. 
Any pseudo-Anosov $\Phi \in \Psi_P$ occurs as the monodromy of a Dehn filling of one of the $M_k$. 
In particular, there exists a constant $V = V(P)$ such that 
$\mathrm{vol}(\Phi) \le V$ holds for any $\Phi \in \Psi_P$. 
\end{thm}

\noindent 
Agol also proved Theorem~\ref{thm_finite-FLM} by using periodic splitting sequences of pseudo-Anosov mapping tori \cite{Agol2}. 
By Theorem~\ref{thm_finite-FLM},  one sees  that 
the following sets $\mathcal{U}$, $\mathcal{U}^+$ and $\mathcal{V}$ are finite.  
\begin{eqnarray*}
&\mathcal{U}&= \{{\Bbb T}(\Phi|_{\varSigma^{\circ}})\ |\ g \ge 2, \ \Phi\ \mbox{is\ a\ pseudo-Anosov\ homeomorphism\ on\ } \varSigma= \varSigma_g \mbox{\ such\ that\ }\lambda(\Phi)= \delta_g\}, 
\\
&\mathcal{U}^+&= \{{\Bbb T}(\Phi|_{\varSigma^{\circ}})\ |\ g \ge 2, \ \Phi\ \mbox{is\ an\ orientable\  pseudo-Anosov\ homeomorphism} 
\\
&\ &  \hspace{2cm}\mbox{on\ }\varSigma= \varSigma_g \mbox{\ such\ that\ }\lambda(\Phi)= \delta_g^+\}, 
\\
&\mathcal{V}&= \{{\Bbb T}(\Phi|_{\varSigma^{\circ}})\ |\ n \ge 3, \ \Phi\ \mbox{is\ a\ pseudo-Anosov\ homeomorphism\ on\ } \varSigma= D_n \mbox{\ such\ that\ } \lambda(\Phi)= \delta(D_n)\}, 
\end{eqnarray*}
where $\delta(D_n)$ denotes the minimal dilatation for pseudo-Anosov elements of  $\mathrm{Mod}(D_n)$ on an 
$n$-punctured disk $D_n$.

The previous study \cite{KT} by the authors implies that $N \in \mathcal{V}$. 
In fact, the mapping class $\phi \in \mathrm{Mod}(D_4)$ represented by the $4$-braid $\sigma_1 \sigma_2 \sigma_3^{-1}$ 
has the minimal dilatation $\delta(D_4)$ \cite{KLS}. 
For the pseudo-Anosov representative $\Phi $ of this mapping class $\phi$, the mapping torus 
${\Bbb T}(\Phi|_{\varSigma^{\circ}})$ is homeomorphic to $N$ \cite{KT}.  
Moreover for each $n \ge 6$ (resp. $n= 3,4,5$), 
a pseudo-Anosov homeomorphism $\Phi_n :D_n \rightarrow D_n$  having the smallest known dilatation (resp. smallest dilatation) 
occurs as the monodromy on a particular fibration on a manifold obtained from $N$ by Dehn filling \cite{KT}. 
See also work of Venzke \cite{Venzke}.  

Hironaka obtained  Theorems~\ref{thm_bound_H}  and \ref{thm_asymptotic} 
by viewing the monodromies of fibrations on manifolds obtained from the $6_2^2$ link exterior 
$S^3 \setminus 6_2^2$ by Dehn filling two cusps. 
(For the link $6_2^2$,  see Figure~\ref{fig_3chain}(right) or Rolfsen's table \cite[Appendix C]{Rolfsen}.)  
There exists an orientable  monodromy $: \varSigma_2 \rightarrow \varSigma_2$ with dilatation $\delta_2= \delta_2^+$ 
of a fibration on a manifold obtained from $S^3 \setminus 6_2^2$ by Dehn filling two cusps. 
This implies that $S^3 \setminus 6_2^2 \in \mathcal{U} \cap  \mathcal{U}^+$ (Lemma~\ref{lem_list-H} or \cite{Hironaka}). 
We see that $S^3 \setminus 6_2^2$ is homeomorphic to  $N(\tfrac{-1}{2})$ (see \cite[Table A.1]{MP} for example), 
where $N(r)$  is the manifold obtained from $N$ by Dehn filling one cusp along the slope $r$. 
As mentioned, computer experiments say that the pseudo-Anosovs  having small dilatation have small volume, and 
$N$ is the candidate having the smallest volume among orientable $3$-manifolds with $3$ cusps. 
These results led us to see monodromies of fibrations on manifold obtained from $N$ by Dehn filling.

In this paper, we investigate the fibrations on manifolds 
obtained from the three $2$-cusped manifolds 
$N(\tfrac{-1}{2})$, $N(\tfrac{-3}{2})$ and $N(2)$ by Dehn filling $2$ cusps. 
The second one $N(\tfrac{-3}{2})$ is homeomorphic to $N(-4)$ 
and this is the Whitehead sister link exterior, i.e, the $(-2,3,8)$-pretzel link exterior (see \cite[Table A.1]{MP}), see Figure~\ref{fig_3chain}(center). 
The manifold $N(\tfrac{-3}{2})$ and the Whitehead link exterior have the smallest volume among orientable $2$-cusped hyperbolic  $3$-manifolds  \cite{Agol}. 
We shall see that 
$N(\tfrac{-3}{2})$ and $N(2)$ are elements of $\mathcal{U}^+$ (Lemmas~\ref{lem_list-KT1}, \ref{lem_list-KT2}).  
Our main result is that $N(\tfrac{-3}{2})$ (resp. $N(2)$) also admits Dehn fillings giving a sequence of fibers over the circle, 
with closed fibers $\varSigma_g$ of genus $g$ for each $g \ge 3$ 
such that the monodromies associated to the fibrations satisfy the same asymptotic inequality as Theorem~\ref{thm_asymptotic}. 
More precisely, we shall prove the following.

\begin{thm}
\label{thm_three}
Let  $r \in \{ \tfrac{3}{-2}, \tfrac{1}{-2}, 2\}$. 
For each  $g \ge 3$,  
there exist  $\varSigma_g$-bundles over the circle 
obtained from  $N(r)$  by Dehn filling all two cusps  
along the boundary slopes of  fibers of  $N(r)$. 
Among them,  
there exist the monodromies  
$ \Phi_g(r): \varSigma_g \rightarrow \varSigma_g$ of the fibrations such that  
\begin{description}
\item[(1)] 
$\displaystyle \lim_{g \to \infty} g  \log \lambda (\Phi_g) = \log( \tfrac{3+ \sqrt{5}}{2})$, 
\item[(2)] 
$\displaystyle \lim_{g \to \infty} \mathrm{vol}(\Phi_g)= \mathrm{vol}(N(r))$.
\end{description}
%where $ \mathrm{vol}(N(r))$ is the hyperbolic volume of  $N(r)$. 
\end{thm} 

\noindent
Independently, Aaber and Dunfield have investigated $\varSigma_g$-bundles over the circle 
obtained from $N(\tfrac{-3}{2})$  by Dehn filling two cusps, see \cite{AD} and  Remark~\ref{rem_AD}. 
 They have obtained similar results on the dilatation to those given in this paper. 
Theorem~\ref{thm_three} in the case $r=  \tfrac{-3}{2}$ was also established by \cite{AD}.

By using monodromies on closed fibers coming from  $N(\tfrac{-3}{2})$, we find an upper bound of $\delta_g$.

\begin{thm} 
\label{thm_bound_KT1}
\begin{description}
\item[(1)] 
$\delta_g \le \lambda_{(g+2,1)} $  if $g \equiv 0,1,5,6 \pmod {10}$ and $g \ge 5$. 
\item[(2)] 
$\delta_g \le \lambda_{(g+2,2)} $ if $g \equiv 7,9 \pmod {10}$ and $g \ge 7$. 
\end{description}
\end{thm}

\begin{thm}
\label{thm_bound_KT2}
Let  $g \equiv 2,4 \pmod{10}$. 
Suppose that $g+2 \not \equiv 0$ $\pmod{4641 (= 3 \cdot 7 \cdot 13 \cdot 17 )}$.  
\begin{description}
\item[(1)] 
$\delta_g \le \lambda_{(g+2,3)} $ if $\gcd(g+2,3)=1$. 
\item[(2)] 
$\delta_g \le \lambda_{(g+2,7)} $ if $3$ divides $g+2$ and $\gcd(g+2,7)=1$. 
\item[(3)] 
$\delta_g \le \lambda_{(g+2,13)} $ if $21 (=3 \cdot 7)$ divides $g+2$ and $\gcd(g+2,13)=1$. 
\item[(4)] 
$\delta_g \le \lambda_{(g+2,17)} $ if $273 (=3 \cdot 7 \cdot 13)$ divides $g+2$ and $\gcd(g+2,17)=1$. 
\end{description}
\end{thm}

 \noindent
 We will verify the bounds in Theorems~\ref{thm_bound_KT1}, \ref{thm_bound_KT2} are sharper than the ones in Theorem~\ref{thm_bound_H} 
 (see  Propositions~\ref{prop_bound_KT1}(1),(2) and \ref{prop_bound_KT2}).
Theorems~\ref{thm_bound_KT1}, \ref{thm_bound_KT2} do not include the case $g \equiv 3,8 \pmod {10}$. 
This is because in this case, $N(\tfrac{-3}{2})$ can not give rise to the monodromy on a closed fiber of genus $g$ whose dilatation is strictly smaller than the one obtained from 
$N(\tfrac{-1}{2})$, see Proposition~\ref{prop_bound_KT1}(3),(4). 
However in  case $g = 8,13$, we find a sharper upper bound than the one in Theorem~\ref{thm_bound_H}. 
Let $\lambda_{(x, y,z)}$ be the largest real root of the polynomial 
$$f_{(x,y,z)}(t)= t^{x+y-z}-t^x - t^y - t^{x-z}- t^{y-z}+1.$$
\begin{prop}
\label{prop_small_genus}
\begin{description}
\item[(1)] 
$\delta_8 \le \lambda_{(18,17,7)} (\thickapprox 1.10403) < \lambda_{(9,1)} (\thickapprox 1.11350)$. 
\item[(2)] 
$\delta_{13} \le \lambda_{(27,21,8)} (\thickapprox 1.07169) < \lambda_{(14,3)} (\thickapprox 1.07266)$. 
\end{description}
\end{prop}

We turn to the study on  $\delta_g^+$. 
We record results by Lanneau-Thiffeault.

\begin{thm}[\cite{LT}] 
\label{thm_lower_bound}
The minimal dilatation $\delta_g^+$ for $ g =6,7$ is not less than the largest real root of the following polynomial. 
\begin{description}
\item[(1)] 
$f_{(6,1)}(t) = t^{12}-t^7-t^6-t^5+1$ if $g= 6$. $(\delta_6^+ \ge \lambda_{(6,1)}  \thickapprox 1.17628.)$ 
\item[(2)] 
$f_{(9,2)}(t)= (t^4-t^3+t^2-t+1)(t^{14}+t^{13}-t^9-t^8-t^7-t^6-t^5+t+1)$ if $g=7$. $(\delta_7^+ \ge \lambda_{(9,2)} \thickapprox 1.11548.)$ 
%\item[(3)] 
%$f_{(8,1)}(t)= t^{16}-t^9-t^8-t^7+1$ if $g= 8$. $(\lambda_{(8,1)} \thickapprox 1.12876.)$ 
\end{description}
\end{thm}

\noindent
%Theorem~\ref{thm_lower_bound}(3) and Theorem~\ref{thm_ori_bound_H}(2) imply that $\delta_8^+= \lambda_{(8,1)}$. 
Lanneau-Thiffeault asked the following.  

\begin{ques}[\cite{LT}]
\label{ques_LT}
For $g$ even,  is $\delta_g^+$ equal to the largest real root of the polynomial 
$$f_{(g,1)}(t)= t^{2g}- t^{g+1}-t^g-t^{g-1}+1?$$
Namely, is  $\delta_g^+$ equal to $\lambda_{(g,1)}$ for $g$ even? 
\end{ques}

\noindent
An upper bound of $\delta_g^+$ given by Hironaka is as follows.  

\begin{thm}[\cite{Hironaka}] 
\label{thm_ori_bound_H}
\begin{description}
\item[(1)] 
$\delta_g^+ \le \lambda_{(g+1,3)}$ if $g \equiv 1,3 \pmod 6$. 
\item[(2)] 
$\delta_g^+ \le \lambda_{(g,1)}$ if $g \equiv 2,4 \pmod 6$. 
\item[(3)] 
$\delta_g^+ \le \lambda_{(g+1,1)}$ if $g \equiv 5 \pmod 6$. 
\end{description}
\end{thm}

\noindent
We do not know whether there exists an orientable pseudo-Anosov homeomorphism of genus $g$ having the dilatation  
$\lambda_{(g,1)}$ (appeared in Question~\ref{ques_LT})  or not for each $g \equiv 0 \pmod 6$. 
Under the assumption that Question~\ref{ques_LT} is true, 
the inequality $\delta_g^+ \le \delta_{g+1}^+$ holds whenever $g \equiv 5 \pmod 6$ and  $\delta_g < \delta_g^+$ holds for all even $g$, see \cite{Hironaka}. 

We give an  upper bound of $\delta_g^+$ in the case $g \equiv 1,5,7,9 \pmod{10}$ using 
orientable pseudo-Anosov monodromies  coming from $N(\tfrac{-3}{2})$.

\begin{thm}
\label{thm_ori_bound_KT}
\begin{description}
\item[(1)] 
$\delta_g^+  \le \lambda_{(g+2,2)}$ 
if $g \equiv 7,9 \pmod{10}$ and $g \ge 7$.  
\item[(2)]  
$\delta_g^+  \le \lambda_{(g+2,4)}$ 
if $g \equiv 1,5 \pmod{10}$ and $g \ge 5$.  
\end{description}
\end{thm}

\noindent
We shall see that 
the bound in Theorem~\ref{thm_ori_bound_KT} improves the one in Theorem~\ref{thm_ori_bound_H} (see Proposition~\ref{prop_ori_bound_KT}). 
Theorem~\ref{thm_ori_bound_KT}(1) together with Theorem~\ref{thm_lower_bound}(2) gives: 

\begin{cor}
\label{cor_genus-7}
$\delta_7^+= \lambda_{(9,2)}$.
\end{cor}

\noindent
Independently, Corollary~\ref{cor_genus-7} was established by Aaber and Dunfiled \cite{AD}.

The following tells us that the sequence  $\{\delta_g^+\}_{g \ge 2}$ is not monotone decreasing.

%which says that $\{\delta_g^+\}_{g \ge 2}$ is not monotone decreasing: 
\begin{prop}
\label{prop_non-monotone}
If Question~\ref{ques_LT} is true, then 
$\delta_{g}^+ < \delta_{g+1}^+$ whenever  $g \equiv 1,5,7,9 \pmod {10}$ and $g \ge 7$. 
In particular the inequality $\delta_7^+ < \delta_8^+$ holds. 
\end{prop}

\noindent
Our pseudo-Anosov homeomorphisms providing  the upper bound of $\delta_g  $ in Theorem~\ref{thm_bound_KT1}(1) are not orientable 
(Remark~\ref{rem_our-example}). 
This together with the inequality $\lambda_{(7,1)} < \lambda_{(6,1)} = \delta_5^+$ implies: 

\begin{cor}
$\delta_5 < \delta_5^+$. 
\end{cor}

We have a question: 

\begin{ques}
\label{ques_magic}
Does the magic manifold $N$ satisfy the following  (1),(2) and (3)? 
\begin{description}
\item[(1)] 
There exist Dehn fillings of $N$ giving an infinite sequence of fiberings over the circle, with closed fibers $\varSigma_{g_i}$ of genus $g_i \ge 2$ 
with $g_i \to \infty$, and with monodromy $\Phi_i$ so that 
$\delta_{g_i} = \lambda(\Phi_i)$. 
\item[(2)] 
There exist Dehn fillings of $N$ giving an infinite sequence of fiberings over the circle, with closed fibers $\varSigma_{g_i}$ of genus $g_i \ge 2$ 
with $g_i \to \infty$, and with monodromy $\Phi_i$ having the orientable (un)stable foliation  so that 
$\delta_{g_i}^+ = \lambda(\Phi_i)$. 
\item[(3)] 
There exist Dehn fillings of $N$ giving an infinite sequence of fiberings over the circle, with fibers $D_{n_i}$ having  $n_i$ punctures 
with  $n_i \to \infty$, and with monodromy $\Phi_i$ so that 
$\delta(D_{n_i}) = \lambda(\Phi_i)$. 
\end{description}
\end{ques}

\noindent
The existence of the manifold satisfying each of  (1),(2) and (3) is guaranteed from Theorem~\ref{thm_finite-FLM}. 
Question~\ref{ques_magic} asks whether $N$ enjoys all (1),(2) and (3)  or not.

The paper is organized as follows. 
We review basic facts in Section~\ref{section_notation}. 
The fibered faces and the entropy function for $N$ are described  in Section~\ref{section_magic-manifold}. 
The (un)stable foliation for the monodromy of the fibration   on  $N$ is discussed in the section. 
We prove theorems in Section~\ref{section_hyperbolic}.

\section{Notation and basic facts} 
\label{section_notation}

\subsection{Pseudo-Anosov}

The {\it mapping class group} $\mathrm{Mod}(\varSigma)$ is the group of isotopy classes of orientation preserving homeomorphisms of an orientable surface $\varSigma$, 
where the group operation is induced by composition of homeomorphisms. 
An element of this group  is called a {\it mapping class}.

A homeomorphism $\Phi: \varSigma \rightarrow \varSigma$ is {\it pseudo-Anosov}  
 if  there exists a constant $\lambda= \lambda(\Phi)>1$ called the {\it dilatation of} $\Phi$  
 and there exists a pair of transverse measured foliations $\mathcal{F}^s$ and $\mathcal{F}^u$ such that 
 $$\Phi(\mathcal{F}^s)= \tfrac{1}{\lambda} \mathcal{F}^s \ \mbox{and}\  \Phi(\mathcal{F}^u)= \lambda \mathcal{F}^u.$$ 
 Measured foliations $\mathcal{F}^s$ and $\mathcal{F}^u$ are called the {\it stable} and {\it unstable foliations} or {\it invariant foliations}  for $\Phi$. 
 In this case the mapping class  $\phi=[\Phi]$ is called pseudo-Anosov. 
 We define the dilatation of $\phi$, denoted by $\lambda(\phi)$, to be the dilatation of $\Phi$.

 The (topological)  entropy $\mathrm{ent}(f)$ is a measure of the complexity of a continuous self-map $f$ on a compact manifold, see for instance  \cite{Walters}. 
 The inequality 
 $$\log \mathrm{sp}(f_{*}) \le \mathrm{ent}(f)$$ 
 holds (see  \cite{Manning}), 
 where  $\mathrm{sp}(f_{*}) $ is the spectral radius of the induced map $f_{*}: H_1(S;{\Bbb R}) \rightarrow H_1(S;{\Bbb R})$ 
 on the first homology group. 
 For any pseudo-Anosov homeomorphism $\Phi: \varSigma \rightarrow \varSigma$, 
the equality 
$$\mathrm{ent}(\Phi)= \log (\lambda(\Phi))$$ 
holds  and $\mathrm{ent}(\Phi)$ attains the minimal entropy among all  homeomorphisms  
which are isotopic to $\Phi$, see \cite[Expos\'{e} 10]{FLP}.  
We denote by  $\mathrm{ent}(\phi)$,  this characteristic number.  
If $\Phi$ has orientable invariant foliations, then the equality 
$$\log \mathrm{sp}(\Phi_{*}) = \mathrm{ent}(\Phi)$$ 
holds, see \cite{Rykken}. 
The converse is  true: 

\begin{thm}[\cite{LT}]
\label{thm_spec}
A pseudo-Anosov homeomorphism $\Phi$ is orientable if and only if 
$\mathrm{sp}(\Phi_{*}) = \lambda(\Phi)$.  
\end{thm}

If we fix a surface $\varSigma$ and take a constant $c>1$, then 
the set of dilatations $\lambda(\Phi)<c$ for pseudo-Anosov homeomorphisms $\Phi: \varSigma \rightarrow \varSigma$ 
is finite, see \cite{Ivanov}. 
In particular the set 
$$\mathrm{Dil}(\varSigma) = \{\lambda(\phi)\ |\ \mbox{pseudo-Anosov\ }\phi \in \mathrm{Mod}(\varSigma)\}$$ 
achieves a minimum $\delta(\varSigma)$.

Thurston's hyperbolization theorem \cite{Thurston3} asserts that 
$\phi$ is pseudo-Anosov if and only if its mapping torus 
$${\Bbb T}(\phi) =\varSigma \times [0,1]/ \sim ,$$ 
 where $\sim$ identifies $(x,1)$ with $(f(x),0)$ for a representative $f \in \phi$, is hyperbolic. 
 We denote the hyperbolic volume of ${\Bbb T}(\phi)$ by $\mathrm{vol}(\phi)$.

Let us suppose that  $\varSigma$ is a compact orientable surface of genus $g$ and we consider
a pseudo-Anosov homeomorphism $\Phi: \varSigma \rightarrow \varSigma$. 
The stable foliation  for $\Phi$ is denote by $\mathcal{F}$.  
Let $x_1, \cdots, x_m$ be all the  singularities for $\mathcal{F}$ in the interior $int(\varSigma)$, and 
$p(x_i) \ge 3$ denotes the number of prongs of $\mathcal{F}$ at $x_i$. 
Let $y_1, \cdots, y_n$ be all the singularities for $\mathcal{F}$ on the boundary $\partial \varSigma$, and 
$p(y_j) \ge 1$ denotes the number of prongs of $\mathcal{F}$ at $y_j$. 
The following Euler-Poincar\'{e} formula holds: 
$$\sum_{i=1}^m (p(x_i)-2) + \sum_{j=1}^n (p(y_j)-2) = -2\chi(\varSigma_g)= 4g-4$$
(see \cite[Expos\'{e} 5]{FLP} for example). 
The pair of integers  
$$(p(x_1)-2, p(x_2)-2, \cdots,  p(x_m)-2,  p(y_1)-2,  p(y_2)-2,\cdots, p(y_n)-2)$$ 
is called the {\it singularity data} of $\Phi$.

\subsection{Thurston norm}

Let $M$ be an irreducible, atoroidal and oriented $3$-manifold with boundary $\partial M$ (possibly $\partial M = \emptyset$). 
We recall the Thurston norm  $\| \cdot \|: H_2(M, \partial M; {\Bbb R}) \rightarrow {\Bbb R}$  (see  \cite{Thurston1}).  
The norm $\| \cdot \| $ has the property that 
for any integral class $a  \in H_2(M, \partial M; {\Bbb R})$, 
$$\| a\|= \min_F \{- \chi(F)\},$$  
where 
%$\chi(F)$ is the  Euler characteristic of $F$ and 
 the minimum is taken over all  oriented surface $F$ embedded in $M$, satisfying $a= [F]$,  
 with no components  of non-negative Euler characteristic. 
The surface $F$ which realizes this minimum  is called  a {\it minimal representative} of $a $.  
For a rational class $a \in H_2(M, \partial M; {\Bbb R})$, take a rational number $r $ so that 
$ra$ is an integral class. 
Then $\|a\|$ is defined to be 
$$\|a\| = \tfrac{1}{|r|} \|r\|.$$ 
The function $\| \cdot \|$ defined on rational classes  admits a unique continuous extension to 
$H_2(M, \partial M; {\Bbb R}) $ which is linear on the ray though the origin. 
The unit ball $U= \{a \in H_2(M, \partial M; {\Bbb R})\ |\ \|a\| \le 1\}$ is a compact, convex polyhedron \cite{Thurston1}.

The following notations are needed to describe how fibrations of $M$ are related to $\| \cdot \|$. 

\begin{itemize}
\item A top dimensional face in the boundary $\partial U$  is denoted by $\Delta$, and its open face is denoted by $int(\Delta)$. 
%\item The cone with the origin $0$ over   $ \Delta$ but the origin $0$ is denoted by $C_{\Delta}$. 
%\item The open cone with the origin  over $\Delta$ is denoted by $int(C_{\Delta})$. 
\item $C_{\Delta}$ is the cone over $\Delta$ with the origin and $int(C_{\Delta})$ is its interior. 
\item The set of integral classes (resp. rational classes) of $int(C_{\Delta})$ is denoted by $int(C_{\Delta}(\Bbb {Z})) $ (resp. $int(C_{\Delta}(\Bbb {Q})) $).  
\end{itemize}

\begin{thm}[\cite{Thurston1}] 
\label{thm_norm-fibration}
Suppose that $M$ is a surface bundle over the circle and let $F$ be a fiber. 
Then there exists a top dimensional  face $\Delta$ satisfying the following. 
\begin{description}
\item[(1)] $[F] \in int(C_{\Delta}(\Bbb Z))$. 
\item[(2)]  
For any $a \in int(C_{\Delta}(\Bbb Z))$, 
a minimal  representative  of $a$ is a  fiber of fibrations on $M$.  
\end{description}
\end{thm}
\noindent
The face $\Delta$  in Theorem~\ref{thm_norm-fibration} is called a {\it fibered face} 
and an integral class $a \in int(C_{\Delta})$ is called a {\it fibered class}.

\subsection{Entropy function}  

Let $M$ be a hyperbolic surface bundle over the circle. 
We fix a fibered face $\Delta$ for $M$. 
The {\it entropy function} 
$ \mathrm{ent}: int(C_{\Delta}) \rightarrow {\Bbb R}$ introduced by Fried  \cite{Fried} is defined as follows. 
The minimal  representative $F_a$ for a primitive fibered class $a \in int(C_{\Delta})$  is connected and is a fiber of fibrations on $M$. 
Let  $\Phi_a : F_a \rightarrow F_a$ be the monodromy. 
Since $M$ is hyperbolic,  $\phi_a = [\Phi_a]$ is  pseudo-Anosov. 
The entropy  $\mathrm{ent}(a) $ and  dilatation $\lambda(a)$ are defined to be the entropy and dilatation of $\phi_a$.
For $r \in {\Bbb Q}$, the entropy $\mathrm{ent}(ra) $ is defined by $ \frac{1} {|r|}  \mathrm{ent}(a)$. 
Fried proved that  $\tfrac{1}{\mathrm{ent}}: int(C_{\Delta}({\Bbb Q})) \rightarrow {\Bbb R}$ 
is concave \cite{Fried}, 
and in particular it  admits a unique continuous extension $\mathrm{ent}: int(C_{\Delta}) \rightarrow {\Bbb R}$.   
Moreover, he proved that the restriction of  $\mathrm{ent}$  to  $int(\Delta)$  is proper, 
that is 
$\mathrm{ent}(a)$  goes to  $ \infty$  as  $a$  goes to 
a  boundary point of  $\partial \Delta$.  
Note that  
$\frac{1}{\mathrm{ent}} : int(C_{\Delta}) \rightarrow {\Bbb R}$  is linear along each ray through the origin 
by definition and cannot be strictly concave for this direction. 
However  Matsumoto and  later McMullen proved that it is strictly concave for other directions.  

\begin{thm}[\cite{Matsumoto,McMullen}] 
\label{thm_strictly-concave}
$\tfrac{1}{\mathrm{ent}}: int(\Delta) \rightarrow {\Bbb R}$ 
is strictly concave. 
\end{thm}

\noindent
By definition of $\mathrm{ent}$, 
$\mathrm{Ent}= \|\cdot\| \mathrm{ent}:  int(C_{\Delta}) \rightarrow {\Bbb R}$ 
is constant on each ray in $int(C_{\Delta}) $ through the origin.  
We call $\mathrm{Ent}(a)$ the {\it normalized entropy} of $a \in int(C_{\Delta})$.  
By Theorem~\ref{thm_strictly-concave} together with the properness of  $\mathrm{ent}$  by Fried, 
$\mathrm{Ent}$  admits the unique minimum 
at the unique ray through the origin.  
In other words,    
if we regard  $\mathrm{Ent}$  as a function defined on $int(\Delta)$, 
then it has the minimum at the unique point in  $int(\Delta)$.

The following question was posed by McMullen. 

\begin{ques}[\cite{McMullen}]
\label{ques_mcmullen}
On which ray in $ int(C_{\Delta})$ does  $\mathrm{Ent}$ attain the minimum? 
Is the minimum attained on a rational class of $int(\Delta)$?
\end{ques}

\noindent
We consider Question~\ref{ques_mcmullen} for  $N(\tfrac{-3}{2})$, $N(\tfrac{-1}{2})$ and $N(2)$, 
see Proposition~\ref{prop_mini-ray_2cusps}.

\section{Magic manifold}
\label{section_magic-manifold}

\subsection{Fibered face and entropy function}
\label{subsection_FiberFace}

The magic manifold $N$ is  a  surface bundle over the circle (\cite{KT} for instance). 
In this section, we recall the entropy function on a fibered face for $N$. 

\begin{figure}[htbp]
\begin{center}
\includegraphics[width=4in]{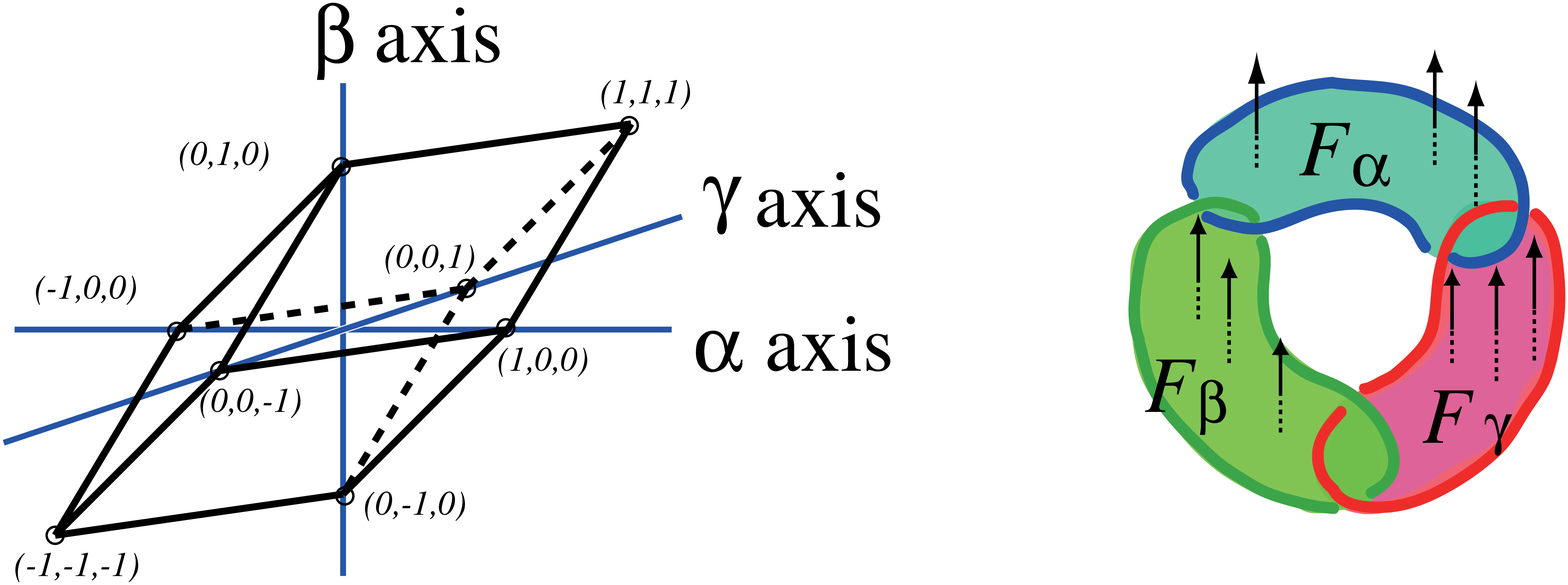}
\caption{(left) Thurston norm ball for N. (right) $F_{\alpha}$, $F_{\beta}$, $F_{\gamma}$. 
[arrows indicate the normal direction of oriented surfaces.]}
\label{fig_poly}
\end{center}
\end{figure}

Let $K_{\alpha}$, $K_{\beta}$ and  $K_{\gamma}$ be the components of  the $3$ chain link $\mathcal{C}_3$. 
They bound the oriented disks $F_{\alpha}$, $F_{\beta}$ and $F_{\gamma}$ with $2$ holes, see Figure~\ref{fig_poly}(right). 
Let $\alpha = [F_{\alpha}]$, $ \beta = [F_{\beta}]$, and $\gamma= [F_{\gamma}]$. 
The Thurston  unit ball $U$ for $N$  is the parallelepiped with vertices 
$\pm \alpha = (\pm 1,0,0)$, $\pm \beta = (0,\pm 1,0)$, $\pm \gamma= (0, 0, \pm1)$, $\pm(\alpha + \beta + \gamma)$, see Figure~\ref{fig_poly}(left). 
The set $\{\alpha, \beta, \gamma\}$ is a basis of $H_2(N,  \partial N; {\Bbb Z})$. 
The symmetry of $\mathcal{C}_3$ tells us  that every top dimensional face is a fibered face. 
We fix a  face $\Delta$ with vertices $\alpha = (1,0,0)$, $\alpha+ \beta + \gamma= (1,1,1)$, $\beta=(0,1,0)$ and $-\gamma= (0,0,-1)$.  
Then 
$$ int(\Delta) = \{x \alpha + y \beta + z \gamma\ |\ x+y-z =1, \ x >0,\  y>0,\   x >z,\   y>z\}.$$ 
Hence if $x \alpha + y \beta + z \gamma  \in int(C_{\Delta})$, then 
\begin{equation}
\label{equation_norm}
\|x \alpha + y \beta + z \gamma\| = x + y -z.
\end{equation}

Let $\mathcal{N}(L)$ be the  regular neighborhood of a link $L$ in $S^3$. 
We denote the tori $ \partial \mathcal{N}(K_{\alpha})$, $ \partial \mathcal{N}(K_{\beta})$, $ \partial \mathcal{N}(K_{\gamma})$ by 
$T_{\alpha}$, $T_{\beta}$, $T_{\gamma}$ respectively. 
Let  $ x \alpha + y \beta + z \gamma $ be a primitive fibered class in $int(C_{\Delta})$. 
The minimal representative of this class is denoted by $F_{x \alpha + y \beta + z \gamma}$ or $F_{(x,y,z)}$. 
Let us put 
$\partial_{\alpha} F_{(x,y,z)} = \partial F_{(x,y,z)} \cap T_{\alpha}$ which consists of  the parallel simple closed curves on $T_{\alpha}$. 
We define $\partial_{\beta} F_{(x,y,z)} $, $\partial_{\gamma} F_{(x,y,z)} \subset  \partial F_{(x,y,z)} $ in the same manner.

\begin{lem}
\label{lem_topological-type}
Let  $ x \alpha + y \beta + z \gamma $ be a primitive fibered class in $int(C_{\Delta})$. 
The number of the boundary components $\sharp(\partial F_{(x,y,z)})$ is equal to  
$\gcd(x,y+z)+ \gcd(y,z+x)+ \gcd(z,x+y)$,  
where $\gcd(0,w)$ is defined by $|w|$.  
More precisely 
\begin{description}
\item[(1)] 
$ \sharp( \partial_{\alpha} F_{(x,y,z)}) = \gcd(x,y+z)$, 
\item[(2)]
$\sharp (\partial_{\beta} F_{(x,y,z)} )=  \gcd(y,z+x)$, 
\item[(3)]
$\sharp( \partial_{\gamma} F_{(x,y,z)}) = \gcd(z,x+y)$. 
\end{description}
\end{lem}

\noindent
{\it Proof.}
We prove  (1). 
The proof of (2),(3) is similar. 
We have the meridian and longitude basis $\{m_{\alpha}, \ell_{\alpha}\}$ for $T_{\alpha}$. 
Similarly we have the bases $\{m_{\beta}, \ell_{\beta}\}$ for $T_{\beta}$ and 
$\{m_{\gamma}, \ell_{\gamma}\}$ for $T_{\gamma}$. 
We consider the long exact sequence of the homology groups of the pair $(N, \partial N)$. 
The boundary map is given by 
\begin{eqnarray*}
\partial_*: H_2(N, \partial N; {\Bbb R})& \rightarrow& H_1(\partial N; {\Bbb R}), 
\\
\alpha &\mapsto& \ell_{\alpha} - m_{\beta}- m_{\gamma}, 
\\
\beta &\mapsto& \ell_{\beta} - m_{\gamma}- m_{\alpha}, 
\\
\gamma &\mapsto& \ell_{\gamma} - m_{\alpha}- m_{\beta}. 
\end{eqnarray*}
Hence 
\begin{equation}
\label{equation_homology}
\partial_*(x \alpha + y \beta + z \gamma) = x \ell_{\alpha} - (y+z) m_{\alpha} + 
y \ell_{\beta} - (z+x) m_{\beta} + z \ell_{\gamma} - (x+y)m_{\gamma}. 
\end{equation}
Since $F_{(x,y,z)}$ is the minimal representative, 
$\partial_{\alpha} F_{(x,y,z)}$ is a union of oriented parallel simple closed curves on $T_{\alpha}$ 
whose homology class  equals  $x \ell_{\alpha} - (y+z) m_{\alpha}  \in H_1(T_{\alpha}; {\Bbb R})$, 
see (\ref{equation_homology}). 
Thus the number of the components of $\partial_{\alpha} F_{(x,y,z)}$ equals $\gcd(x,y+z)$. 
This completes the proof. 
$\Box$
\medskip

\noindent
From the proof of Lemma~\ref{lem_topological-type}, one sees that 
the boundary slope of each simple closed curve of $\partial_{\alpha} F_{(x,y,z)} $ equals $\tfrac{-(y+z)}{x}$. 
Similarly the boundary slope of each component of $\partial_{\beta} F_{(x,y,z)} $ (resp. $\partial_{\gamma} F_{(x,y,z)} $) is given by 
$ \tfrac{-(z+x)}{y}$ (resp. $\tfrac{-(x+y)}{z}$). 
Let us define 
\begin{equation}
\label{equation_slopeXYZ}
\mathrm{slope}(x \alpha + y \beta + z \gamma) = (\tfrac{-(y+z)}{x}, \tfrac{-(z+x)}{y}, \tfrac{-(x+y)}{z}).
\end{equation}
This notation $\mathrm{slope}(\cdot)$ is needed for the study of Dehn fillings of $N$ in Section~\ref{section_hyperbolic}.

One can compute the entropy for any element of  $ int(C_{\Delta}(\Bbb Z))$ by using the next theorem. 

\begin{thm}[\cite{KT}]
\label{thm_Fried-Oertel-Poly}
The dilatation $\lambda_{(x,y.z)}$ of $x \alpha + y \beta + z \gamma \in int(C_{\Delta}(\Bbb Z))$ is 
the largest real root of the polynomial 
$$P(t_1,t_2,t_3)= -t_1 -t_2 + t_3 + t_1 t_2 -t_1 t_3 -t_2 t_3.$$
\end{thm}

\noindent
Since $P(t^x, t^y, t^z)= t^z(t^{x+y-z}-t^x - t^y - t^{x-z}- t^{y-z}+1)$, 
$\lambda_{(x,y,z)}$ is the largest real root  of 
$$f_{(x,y,z)}(t)= t^{x+y-z}-t^x - t^y - t^{x-z}- t^{y-z}+1.$$ 

The minimum of $\mathrm{Ent}: int(C_{\Delta}) \rightarrow {\Bbb R}$ is equals to $2 \log(2+ \sqrt{3})$ 
and it is is attained by  $\alpha+ \beta$ \cite{KT}. 
Since  the Thurston norm ball of $N$ has a symmetry, $\min \mathrm{Ent}$ does not depend on fibered faces of $N$.

\subsection{Invariant foliations}

Let $\Phi_{(x,y,z)}$ be the monodromy of the fibration of $N$ associated to a primitive fibered class $x \alpha+ y \beta + z \gamma \in int(C_{\Delta})$ and 
let $F_{(x,y,z)}$ be its fiber. 
We denote  the stable foliation for $\Phi_{(x,y,z)}$ by $\mathcal{F}_{(x,y,z)}$. 
We shall compute the number of prongs at the singularities of  $\mathcal{F}_{(x,y,z)}$. 

\begin{prop}
\label{prop_sing-data}
%Let $x \alpha+ y \beta + z \gamma$ be an element  of  $int(C_{\Delta})(\Bbb Z)$. 
The singularity data of $\Phi_{(x,y,z)}$ is given by 
$$(\underbrace{\tfrac{x}{\gcd(x,y+z)} -2, \cdots, \tfrac{x}{\gcd(x,y+z)} -2}_{\gcd(x,y+z)},  
\underbrace{\tfrac{y}{\gcd(y,x+z)} -2, \cdots, \tfrac{y}{\gcd(y,x+z)} -2}_{\gcd(y,x+z)},   
\underbrace{\tfrac{x+y-2z}{\gcd(z,x+y)} -2, \cdots, \tfrac{x+y-2z}{\gcd(z,x+y)} -2}_{\gcd(z,x+y)}).$$
More precisely $\mathcal{F}_{(x,y,z)}$  has 
\begin{description}
\item[(1)]  
$\tfrac{x}{\gcd(x,y+z)}$ prongs  at each component  of $\partial_{\alpha} F_{(x,y,z)}$, 
\item[(2)] 
$\tfrac{y}{\gcd(y,x+z)}$ prongs  at each component  of  $\partial_{\beta} F_{(x,y,z)}$, 
\item[(3)] 
$\tfrac{x+y-2z}{\gcd(z,x+y)}$ prongs  at each component of $\partial_{\gamma} F_{(x,y,z)}$, and 
\item[(4)] 
no singularities in the interior of $F_{(x,y,z)} $. 
\end{description}
\end{prop}

\noindent
Here we recall the formula of the intersection number $i([c],[c'])$ 
between isotopy classes  of essential simples closed curves $c,c'$ on a torus $T$. 
Let $\tfrac{p}{q}, \tfrac{r}{s}$ be rational numbers or $\tfrac{1}{0}$ with irreducible forms and 
suppose that $\tfrac{p}{q}, \tfrac{r}{s}$ are slopes on $T$ which represent isotopy classes $[c],[c']$ respectively. 
Then $i([c],[c'])= |ps-qr|$. 
\medskip

\begin{figure}[htbp]
\begin{center}
\includegraphics[width=3.5in]{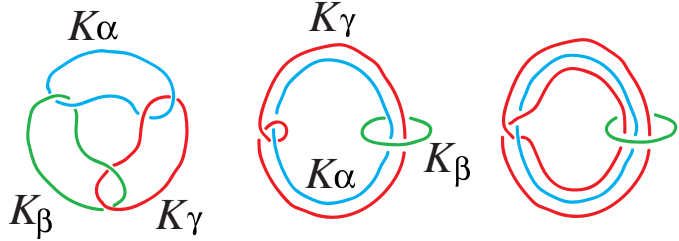}
\caption{(left, center) $\mathcal{C}_3$. (right) $\overline{b}$. 
(this figure explains how to obtain $H$.)}
\label{fig_chain-homeo}
\end{center}
\end{figure}

\begin{figure}[htbp]
\begin{center}
\includegraphics[width=4.5in]{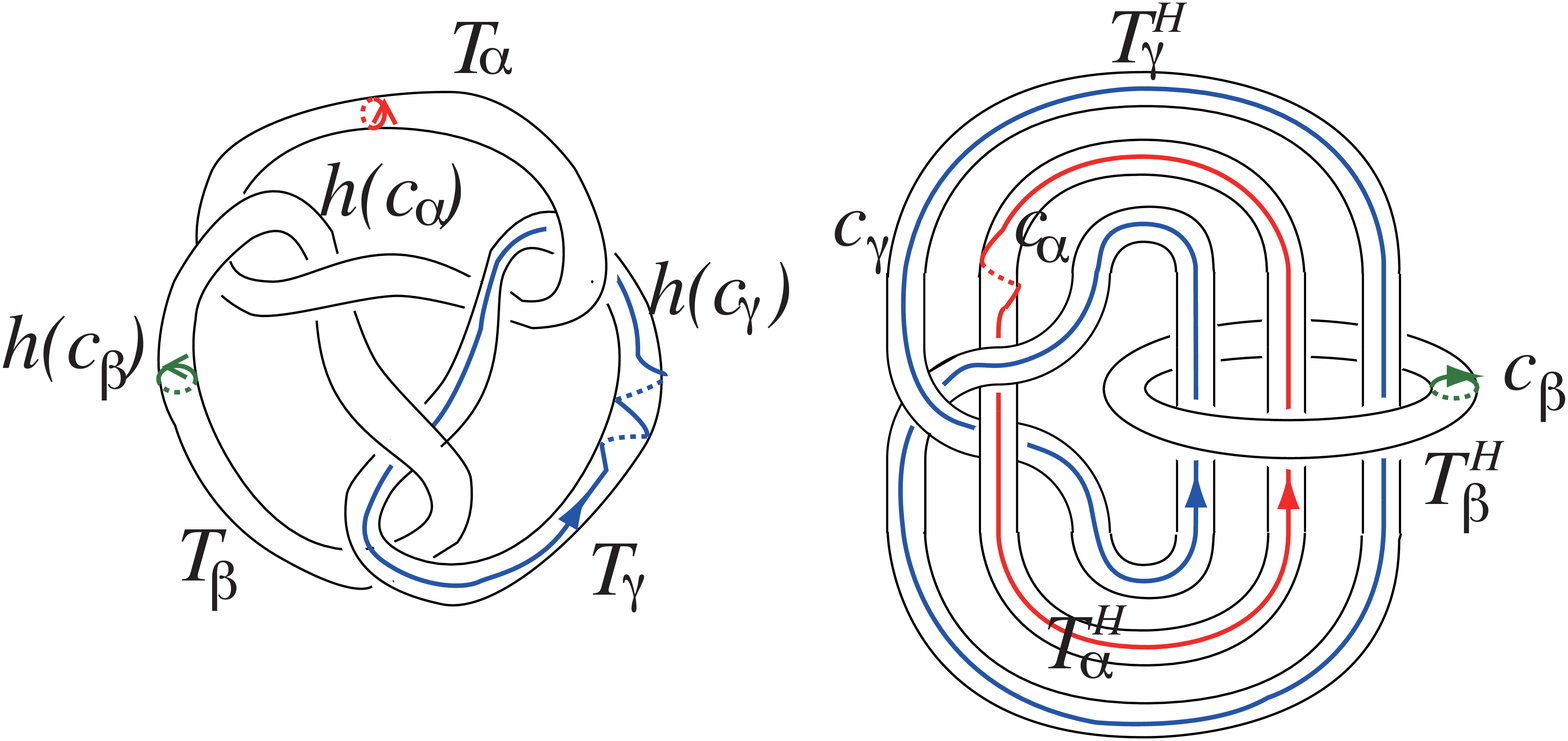}
\caption{(left) $h(c_{\alpha})$, $h(c_{\beta})$, $h(c_{\gamma})$. 
(right) $c_{\alpha}$, $c_{\beta}$, $c_{\gamma}$.}
\label{fig_3braidTori}
\end{center}
\end{figure}

\noindent
{\it Proof of Proposition~\ref{prop_sing-data}.} 
Observe that a fiber $F= F_{(1,1,0)}$ associated to the fibered class $\alpha + \beta$ is  a sphere with $4$ boundary components. 
The monodromy $\Phi= \Phi_{(1,1,0)}: F \rightarrow F$ of the fibration on $N$ is represented by the $3$-braid 
$b= \sigma_2 \sigma_1^{-1} \sigma_2$. 
In particular $S^3 \setminus \overline{b}$ is homeomorphic to $S^3 \setminus \mathcal{C}_3= N$, 
where  $\overline{b}$ is a union of the closed braid of $b$ and the braid axis, see Figure~\ref{fig_chain-homeo}(right). 
We define a homeomorphism  $H: S^3 \setminus \mathcal{N}(\mathcal{C}_3) \rightarrow S^3 \setminus \mathcal{N}(\overline{b})$  as follows. 
Notice that the the link illustrated in Figure~\ref{fig_chain-homeo}(center) is isotopic to $\mathcal{C}_3$. 
 We cut the twice-punctured disk $F_{\alpha}$ bounded by the component $K_{\alpha}$. 
Let $F'_{\alpha}$ and $F''_{\alpha}$ be the resulting twice-punctured disks after cutting $F_{\alpha}$. 
Reglue $F'_{\alpha}$ and $F''_{\alpha}$ after 
twisting the neighborhood of $F'_{\alpha}$ by $360 $ degrees in the clockwise direction. 
Then we obtain the link $\overline{b}$  whose exterior $S^3 \setminus \overline{b}$  is homeomorphic to  $S^3 \setminus \mathcal{C}_3$, 
see Figure~\ref{fig_chain-homeo}. 
The inverse $H^{-1}$ is denoted by $h$. 
We set $T^H_{\alpha}= H(T_{\alpha})$, $T^H_{\beta}= H(T_{\beta})$ and $T^H_{\gamma}= H(T_{\gamma})$, see Figure~\ref{fig_3braidTori}. 
(Then $\partial \mathcal{N}(\overline{b}) = T^H_{\alpha} \cup T^H_{\beta} \cup T^H_{\gamma}$.) 

\begin{figure}[htbp]
\begin{center}
\includegraphics[width=4.5in]{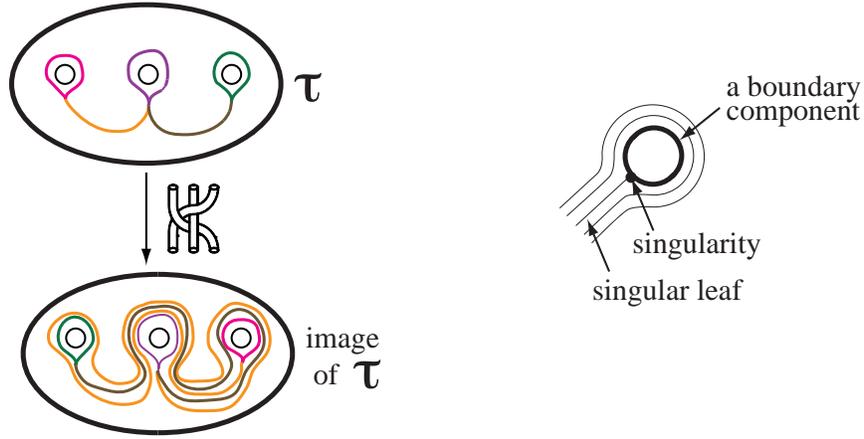}
\caption{(left) invariant train track $\tau$ for $\Phi_{(1,1,0)}$. 
(right) $1$-pronged singularity.}
\label{fig_TrackImage}
\end{center}
\end{figure}

The invariant train track $\tau$ which carries  the stable lamination $\ell^s$ for $\Phi$ is illustrated in Figure~\ref{fig_TrackImage}(left). 
The stable foliation $\mathcal{F}$ for $\Phi$ has a $1$ prong on each component of $\partial F$ and 
it has no singularity in the interior of $F$. 
We consider the suspension flow induced on the mapping torus 
$N=F \times [0,1]/ \sim ,$ 
 where $\sim$ identifies $(x,1)$ with $(\Phi(x),0)$. 
 One obtains the simple closed curve $c_{\alpha} \subset T_{\alpha}^H$  which is the closed orbit of the singularity of $\mathcal{F}$ 
 on $\partial F  \cap T_{\alpha}^H$.  
 Similarly one has the closed orbits $c_{\beta} \subset T_{\beta}^H$, $c_{\gamma} \subset T_{\gamma}^H$, see Figure~\ref{fig_3braidTori}(right). 
 (One can draw these closed orbits by using the orbit of each cusp of $F \setminus \tau$.) 
 Let $\mathcal{L}^s \subset N$ be the suspended stable lamination  which is constructed from $\ell^s \times I \subset F \times I$ 
 by gluing  $\ell^s \times \{1\}$ to $\ell^s \times \{0\}$ using $\Phi$. 
 By construction, $\mathcal{L}^s$ is carried by the branched surface $B_{\tau}$ which is obtained from $\tau \times I$ by gluing 
 $\tau \times \{1\}$ to $\tau \times \{0\}$ using $\Phi$. 
 Notice that  $c_{\alpha}$, $c_{\beta}$ and $c_{\gamma}$ correspond to the branched loci of $B_{\tau}$. 
 By work of Fried \cite{Fried} (see also work of Long-Oertel \cite{LO}), we may assume that  the fiber $F_{(x,y,z)}$ is transverse to $\mathcal{L}^s$. 
 The stable lamination $\ell^s_{(x,y,z)}$ for $\Phi_{(x,y,z)}$ is given by the intersection $\mathcal{L}^s \cap F_{(x,y,z)}$ and 
 $\ell^s_{(x,y,z)}$ is carried by the train track $B_{\tau} \cap F_{(x,y,z)}$. 
 This implies that $\mathcal{F}_{(x,y,z)}$  has no singularity in the interior of $F_{(x,y,z)}$ and 
 we finish the proof of (4). 
 
 We consider the number of prongs of $\mathcal{F}_{(x,y,z)}$ at each component of $\partial_{\alpha} F_{(x,y,z)}$. 
 The boundary slope of each simple closed curve  of $\partial_{\alpha} F_{(x,y,z)}$ is given by $\tfrac{-(y+z)}{x}$. 
 The desired number  is equal to the intersection number 
 $$i([c_{-(y+z)/x}], [h(c_{\alpha})])=  i([H(c_{-(y+z)/x})], [c_{\alpha}]),$$ 
 where $c_r$ is a simple closed curve with slope $r \in {\Bbb Q} \cup \{\tfrac{1}{0}\}$ on $T_{\alpha}$. 
 Observe that $h(c_{\alpha})$ has the slope $\tfrac{1}{0}$ (see Figure~\ref{fig_3braidTori}). 
 Hence 
 $$i([c_{-(y+z)/x}], [h(c_{\alpha})])=|1 \cdot  \tfrac{ x}{\gcd (x,y+z)} +  0 \cdot \tfrac{y+z}{\gcd (x,y+z) }|= \tfrac{x}{\gcd(x,y+z)}.$$
 This completes the proof of (1). 
 
 One verifies that $h(c_{\beta})$ and $h(c_{\gamma})$ have slopes $\tfrac{1}{0}$ and $\tfrac{-2}{1}$ respectively. 
 By using the similar argument, one can prove (2),(3). 
 $\Box$
 \medskip

We consider the orientability of $\mathcal{F}_{(x,y,z)}$ using Theorem~\ref{thm_spec}. 
Alexander polynomial of $\mathcal{C}_3$ is 
$$A(t_1, t_2, t_3)= t_1t_2 + t_2 t_3 + t_3 t_1 - t_1- t_2 - t_3.$$ 
%see \cite[Appendix C]{Rolfsen}. 
The following is a consequence of  Proposition~7.3.10 in \cite{Kawauchi} which tells us that  the relation between the Alexander polynomial of links and 
the characteristic polynomial of  $\Phi_*: H_1(\varSigma;{\Bbb R}) \rightarrow H_1(\varSigma;{\Bbb R})$  on fibers $\varSigma$ in the link exteriors.

\begin{lem}
The spectral radius of $(\Phi_{(x,y,z)})_*$ is the largest absolute value among roots of 
$$A(t^x, t^y, t^z)= t^{x+y}+ t^{y+z}+ t^{z+x}- t^x -t^y - t^z.$$
\end{lem}

\begin{prop}
\label{prop_magic_ori}
The pseudo-Anosov homeomorphism  $\Phi_{(x,y,z)}$ is orientable  if and only if $x$ and $y$ are even and $z$ is odd. 
\end{prop}

\noindent
{\it Proof.} 
(If part.) 
Suppose that $x$ and $y$ are even and $z$ is odd.  
Then 
$$P(t^x, t^y, t^z) = A ((-t)^x, (-t)^y, (-t)^z).$$
This implies that $\lambda(\Phi_{(x,y,z)})=\mathrm{sp}((\Phi_{(x,y,z)})_*)$. 
By Theorem~\ref{thm_spec} $\mathcal{F}_{(x,y,z)}$ is orientable. 
\medskip
\\
(Only if part.) 
Suppose that $x$ or $y$ is odd. 
We may assume that $x$ is odd. 
The number of prongs of $\mathcal{F}_{(x,y,z)}$ at each component of  $\partial_{\alpha} F_{(x,y,z)}$ equals 
$\tfrac{x}{\gcd(x,y+z)}$ which is odd.  
Thus  $\mathcal{F}_{(x,y,z)}$ can not be orientable. 
$\Box$

\subsection{Non-hyperbolic Dehn fillings} 

Let $M$ be a $3$-manifold with boundary tori $T_0, \cdots, T_j$ and 
let $r_i \in {\Bbb Q} \cup \infty $ be a slope on $T_i$. 
Then  $M(r_0, r_1, \cdots, r_j)$ denotes the manifold obtained from $M$ by Dehn filling 
along the slope $r_i$ for each $i$, that is 
$M(r_0, r_1, \cdots, r_j)$ is the manifold  attaching a solid torus $\widetilde{T}_i$ to $M$ along $T_i$ 
in such a way that $r_i$ bounds a disk in $\widetilde{T}_i$.

Martelli and Petronio classified all the non-hyperbolic fillings of the magic manifold \cite[Theorems 1.1, 1.2, 1.3]{MP}. 
We denote by $T_0, T_1, T_2$, the boundary tori of  $N= S^3 \setminus \mathcal{N}(\mathcal{C}_3)$.

\begin{thm}[\cite{MP}]
\label{thm_magic-list}
\begin{description}
\item[(1)] 
$N(\tfrac{p}{q})$ is hyperbolic if and only if 
$$\tfrac{p}{q} \notin \{\infty, -3, -2, -1, 0\}.$$
\item[(2)]  
$N(\tfrac{p}{q}, \tfrac{r}{s})$ is hyperbolic if and only if 
$$\tfrac{p}{q}, \tfrac{r}{s} \notin \{\infty, -3, -2, -1, 0\}\  \mbox{and}\  
(\tfrac{p}{q}, \tfrac{r}{s}) \notin \{(1,1), (-4,\tfrac{-1}{2}), (\tfrac{-3}{2}, \tfrac{-5}{2})\}.$$
\end{description}
\end{thm}

\noindent
As a corollary of Theorem~\ref{thm_magic-list} one has: 

\begin{cor}
\label{cor_magic-list}
If  $N(\tfrac{p}{q}, \tfrac{r}{s}, \tfrac{t}{u})$ is hyperbolic, then 
$$\tfrac{p}{q}, \tfrac{r}{s} \notin \{\infty, -3, -2, -1, 0\}\  \mbox{and}\  
(\tfrac{p}{q}, \tfrac{r}{s}) \notin \{(1,1), (-4,\tfrac{-1}{2}), (\tfrac{-3}{2}, \tfrac{-5}{2})\}.$$
\end{cor}

Let us consider the monodromy $\Phi_{(x,y,z)}: F_{(x,y,z)} \rightarrow F_{(x,y,z)} $ of the fibration on $N$ 
associated  to a primitive fibered class $x \alpha + y \beta + z \gamma \in int(C_{\Delta})$. 
Recall that $\mathrm{slope}(x \alpha + y \beta + z \gamma) = (\tfrac{-(y+z)}{x}, \tfrac{-(z+x)}{y}, \tfrac{-(x+y)}{z})$, see (\ref{equation_slopeXYZ}). 
By capping each boundary component of $F_{(x,y,z)}$, $\Phi_{(x,y,z)}$ extends to the monodromy $\overline{\Phi}_{(x,y,z)}$ 
with a closed fiber $\overline{F}_{(x,y,z)}$ of the fibration on $N(\frac{-(y+z)}{x}, \frac{-(z+x)}{y}, \frac{-(x+y)}{z})$. 
If the stable foliation $\mathcal{F}_{(x,y,z)}$ has no $1$ prong on each component of $\partial F_{(x,y,z)}$, 
then $\overline{\Phi}_{(x,y,z)}$ is pseudo-Anosov and 
$\lambda(\overline{\Phi}_{(x,y,z)})=\lambda(\Phi_{(x,y,z)})$. 
If $\mathcal{F}_{(x,y,z)}$ has a $1$ prong on a component of $\partial F_{(x,y,z)}$, then 
$\overline{\Phi}_{(x,y,z)}$ may not be pseudo-Anosov.
%(Hence $N(\frac{-(y+z)}{x}, \frac{-(z+x)}{y}, \frac{-(x+y)}{z})$ may not be hyperbolic.)  

\section{Hyperbolic Dehn fillings $N(\tfrac{-3}{2})$, $N(\tfrac{-1}{2})$ and $N(2)$} 
\label{section_hyperbolic}

%We investigate the monodromies of the fibrations on $ N(\tfrac{-3}{2})$, $ N(\tfrac{-1}{2})$ and $N(2)$. 
%%The entropy functions for these manifolds have a common property (Lemma~\ref{lem_dilatation_kl}). 
%We shall show that $\min \mathrm{Ent}$ for these manifolds are the same (Proposition~\ref{prop_mini-ray_2cusps}). 

\subsection{Thurston norm balls of $N(\tfrac{-3}{2})$, $N(\tfrac{-1}{2})$ and $N(2)$}
\label{subsection_Thurston}

Let  $N(r)$  be the manifold obtained from  $N$  by Dehn filling 
the cusp specified by $T_{\beta}$  along the slope  $r \in \mathbb{Q}$.  
Then there exists a natural injection  
$\iota_{\beta}: H_2(N(r), \partial N(r)) \rightarrow H_2(N, \partial N)$  
whose image equals 
$$S_{\beta}(r)= \{(x,y,z) \in H_2(N, \partial N)\ |\ -ry= z+x\},$$
see \cite{KKT2}.  
By Theorem~\ref{thm_magic-list}(1), 
$N(r)$  is hyperbolic if and only if  $r \in {\Bbb Q} \setminus  \{-3, -2, -1, 0\}$.  
Choose $r \in {\Bbb Q} \setminus  \{-3, -2, -1, 0\}$,  
and assume that  $a \in S_{\beta}(r) = \mathrm{Im} \, \iota_{\beta}$  is 
a fibered class of  $H_2(N, \partial N)$.  
Then $\overline{a} = \iota_{\beta}^{-1}(a) \in H_2(N(r), \partial N(r))$  is also 
a fibered class of  $N(r)$. 

Similarly, when  $N(r)$  is the manifold obtained from $N$ by Dehn filling the cusp specified 
by  $T_{\alpha}$  or  $T_{\gamma}$  along the slope  $r$, 
one has natural injections, 
\begin{eqnarray*}
	\iota_{\alpha} & : H_2(N(r), \partial N(r)) \to H_2(N, \partial N), \\ 
	\iota_{\gamma} & : H_2(N(r), \partial N(r)) \to H_2(N, \partial N) 
\end{eqnarray*}
such that their images are 
\begin{eqnarray*}
	S_{\alpha}(r) & = \{ (x, y, z) \in H_2(N, \partial N) \, | \, -rx = y+z \}, \\
	S_{\gamma}(r) & = \{ (x, y, z) \in H_2(N, \partial N) \, | \, -rz = x+y \}.  
\end{eqnarray*}

We set 
\begin{eqnarray*}
 \mathfrak{a}&=& 2 \alpha + 2 \beta + \gamma,\hspace{2mm} \mathfrak{b}= \alpha + 2 \beta + 2 \gamma,
 \\
 \mathfrak{p}&=& \alpha + 2 \beta,\ \hspace{9mm}\mathfrak{q}= 2 \beta + \gamma,
 \\
 \mathfrak{r}&=& \alpha + \beta - \gamma,\ \hspace{4mm}\mathfrak{s}= \alpha - \beta.
\end{eqnarray*}
For $k, \ell \in {\Bbb Z}$, we have 
\begin{equation}
\label{equation_slope}
\begin{split}
\mathrm{slope}(k \mathfrak{a}+ \ell \mathfrak{b}) &= (\tfrac{-3k-4\ell}{2k+\ell}, \tfrac{-3}{2}, \tfrac{-4k-3\ell}{k+2\ell}), 
\\
\mathrm{slope}(k \mathfrak{p}+ \ell \mathfrak{q}) &= (\tfrac{-2k-3\ell}{k}, \tfrac{-1}{2}, \tfrac{-3k-2\ell}{\ell}), 
\\
\mathrm{slope}(k \mathfrak{r}+ \ell \mathfrak{s}) &= (\tfrac{\ell}{k+\ell}, \tfrac{-\ell}{k-\ell}, 2),\ \mbox{and}
\end{split}
\end{equation}
\begin{equation}
\label{equation_3families}
\begin{split}
 k \mathfrak{a}+ \ell \mathfrak{b}&= (2k+\ell) \alpha + (2k+ 2\ell) \beta + (k+ 2 \ell) \gamma \in   \iota_{\beta}(H_2(N(\tfrac{-3}{2}), \partial N(\tfrac{-3}{2}))), 
 \\
 k \mathfrak{p} + \ell \mathfrak{q} &= k \alpha + (2k+ 2\ell) \beta + \ell \gamma \in   \iota_{\beta}(H_2(N(\tfrac{-1}{2}), \partial N(\tfrac{-1}{2}))),
 \\
 k \mathfrak{r} + \ell \mathfrak{s} &= (k+ \ell) \alpha + (k-\ell) \beta  - k \gamma \in   \iota_{\gamma}(H_2(N(2), \partial N(2))). 
 \end{split}
\end{equation}
It is easy to check that 
$\{ \overline{\mathfrak{a}}, \overline{\mathfrak{b}}\}$, $\{ \overline{\mathfrak{p}},  \overline{\mathfrak{q}}\}$ and $\{ \overline{\mathfrak{r}},  \overline{\mathfrak{s}}\}$ 
are bases of $H_2(N(\tfrac{-3}{2}), \partial N(\tfrac{-3}{2}); {\Bbb Z})$, $H_2(N(\tfrac{-1}{2}), \partial N(\tfrac{-1}{2}); {\Bbb Z})$ and  $H_2(N(2), \partial N(2); {\Bbb Z})$ 
respectively. 
%\begin{itemize}
%\item 
%$\{ \overline{\mathfrak{a}}, \overline{\mathfrak{b}}\}$ is a basis of $H_2(N(\tfrac{-3}{2}), \partial N(\tfrac{-3}{2}); {\Bbb Z})$, 
%\item
%$\{ \overline{\mathfrak{p}},  \overline{\mathfrak{q}}\}$ is a basis of $H_2(N(\tfrac{-1}{2}), \partial N(\tfrac{-1}{2}); {\Bbb Z})$, 
%\item 
%$\{ \overline{\mathfrak{r}},  \overline{\mathfrak{s}}\}$ is a basis of $H_2(N(2), \partial N(2); {\Bbb Z})$. 
%\end{itemize}

Note that 
$\gcd(k,\ell)=1$ if and only if 
$ k \mathfrak{a} + \ell \mathfrak{b}$, $ k \mathfrak{p} + \ell \mathfrak{q}$ and  $ k \mathfrak{r} + \ell \mathfrak{s}$ are primitive integral classes of $H_2(N, \partial N; {\Bbb R})$. 
All $k \mathfrak{a}+ \ell \mathfrak{b}, k \mathfrak{p} + \ell \mathfrak{q}, k \mathfrak{r} + \ell \mathfrak{s} $ are fibered classes in  $  int(C_{\Delta}) $ 
for  $k >0$ and $-k < \ell< k$. 

We first focus on the topological types of fibers for primitive fibered classes in $ int(C_{\Delta})$. 
Let $ \varSigma_{g,p}$ be a compact orientable surface of genus $g$ with $p$ boundary components.

\begin{lem}
\label{lem_type_fiber}
Suppose that $k >0$, $-k < \ell< k$ and $\gcd(k,\ell)= 1$. 
\begin{description}
\item[(1)] 
$F_{ k \mathfrak{a}+ \ell \mathfrak{b}}= \varSigma_{k-2, k+ \ell+6}$  if $\gcd (2k+\ell,5)=5$ or $\gcd (5,k+2\ell)=5$. 
Otherwise $F_{ k \mathfrak{a}+ \ell \mathfrak{b}}= \varSigma_{k,k+\ell+2}$. 
\item[(2)] 
$F_{k \mathfrak{p} + \ell \mathfrak{q} }= \varSigma_{k-1, k+ \ell+4}   $ if $\gcd(k,3)=3$ or $\gcd(3,\ell)=3$. 
Otherwise $F_{k \mathfrak{p} + \ell \mathfrak{q}}= \varSigma_{k,k+\ell+2}  $. 
\item[(3)] 
$F_{k \mathfrak{r} + \ell \mathfrak{s} }=\varSigma_{k,k+2}$. 
\end{description}
\end{lem}

\noindent
{\it Proof  of (1).}  
By Lemma~\ref{lem_topological-type}, 
\begin{eqnarray*}
\sharp(\partial F_{ k \mathfrak{a}+ \ell \mathfrak{b}})
&=& \gcd(2k+\ell, 3k+4\ell) + \gcd (2k+ 2\ell, 3k + 3 \ell)+ \gcd(4k+3 \ell, k + 2\ell)
\\
&=& \gcd(2k+\ell, 5k) + k+ \ell+ \gcd(5 \ell, k + 2\ell)
\\
&=& \gcd(2k+\ell, 5) + k+ \ell+ \gcd(5, k + 2\ell). 
\end{eqnarray*}
The last equality holds since  $\gcd(k,\ell)= 1$. 
The following 3 cases can occur. 
\begin{description}
 \item[(1)] 
  $\gcd (2k+\ell,5)=1$ and $\gcd(5,k+2\ell)=1$. 
\item[(2)] 
 $\gcd (2k+\ell,5)=5$ and $\gcd(5,k+2\ell)=1$. 
\item[(3)] 
 $\gcd (2k+\ell,5)=1$ and $\gcd (5,k+2\ell)=5$. 
\end{description}
In the case (1), the genus $g$ of $F_{ k \mathfrak{a}+ \ell \mathfrak{b}}$ must satisfy 
$$-(2-2g-k-\ell-2) = \| k \mathfrak{a}+ \ell \mathfrak{b} \|= 3k + \ell$$
(see (\ref{equation_norm})). 
Thus $g= k$ and  $F_{ k \mathfrak{a}+ \ell \mathfrak{b}}= \varSigma_{k,k+\ell+2}$. 
In the cases (2) and (3), $F_{ k \mathfrak{a}+ \ell \mathfrak{b}}= \varSigma_{k-2, k+ \ell+6}$. 

The proof of claims (2), (3) of the lemma is similar. 
$\Box$

\begin{lem}
\label{lem_symmetry}
Suppose that $k >0$ and $- k < \ell < k$. 
\begin{description}
\item[(1)] 
$F_{ k \mathfrak{a}+ \ell \mathfrak{b}} = F_{ \ell \mathfrak{a}+ k \mathfrak{b}}$ and 
$\lambda(k \mathfrak{a}+ \ell \mathfrak{b}) = \lambda(\ell \mathfrak{a}+ k \mathfrak{b})$. 
\item[(2)] 
$F_{ k \mathfrak{p}+ \ell \mathfrak{q}} = F_{ \ell \mathfrak{p}+ k \mathfrak{q}}$ and 
$\lambda(k \mathfrak{p}+ \ell \mathfrak{q}) = \lambda(\ell \mathfrak{p}+ k \mathfrak{q})$. 
\end{description}
\end{lem}

\noindent
{\it Proof.} 
(1) By the symmetry of the Thurston norm ball of $N$, it is not hard to see that 
the topological type of the minimal representative (resp. the dilatation) for 
$\ell \mathfrak{a}+ k \mathfrak{b} = (2 \ell +k) \alpha + (2 \ell + 2 k) \beta + (\ell + 2k) \gamma$ 
is same as the one for $k \mathfrak{a}+ \ell \mathfrak{b} = (2k+\ell) \alpha + (2 k + 2 \ell) \beta + (k+ 2\ell) \gamma$. 

The proof of (2) is similar. 
$\Box$
\medskip

\noindent
We make a remark that 
it is not true in general that $F_{ k \mathfrak{r}+ \ell \mathfrak{s}} =  F_{ \ell \mathfrak{r}+ k \mathfrak{s}}$ and 
$\lambda(k \mathfrak{r}+ \ell \mathfrak{s}) = \lambda(\ell \mathfrak{r}+ k \mathfrak{s})$ 
for $k >0$ and $- k < \ell < k$. 
We do not use this remark in the rest of the paper.

\begin{lem}
\label{lem_symmetry-fiber}
Suppose that $0 < \ell< k$ and $\gcd(k,\ell)= 1$. 
\begin{description}
\item[(1)] 
The genus of $F_{ k \mathfrak{a}+ \ell \mathfrak{b}}$ equals the one of $ F_{ k \mathfrak{a}- \ell \mathfrak{b}}$. 
\item[(2)] 
The genus of $F_{ k \mathfrak{p}+ \ell \mathfrak{q}}$ equals the one of $ F_{ k \mathfrak{p}- \ell \mathfrak{q}}$. 
\item[(3)] 
The genera of $F_{ k \mathfrak{r}+ \ell \mathfrak{s}}$ and $ F_{ k \mathfrak{r}- \ell \mathfrak{s}}$ equal $k$
\end{description}
\end{lem}

\noindent
{\it Proof.} 
(1) 
By Lemma~\ref{lem_type_fiber}(1), 
the genus of $F_{ k \mathfrak{a} \pm  \ell \mathfrak{b}}$  equals $k-2$ if 
$\gcd(2k \pm \ell, 5)=5$ or $\gcd(5, k \pm  2 \ell)=5$. 
Otherwise its genus equals $k$. 
It is easy to check that 
\begin{itemize}
\item 
$\gcd(2k+ \ell, 5)=5$ if and only if $\gcd(5, k - 2\ell)=5$. 
\item 
$\gcd(2k- \ell, 5)=5$ if and only if $\gcd(5, k + 2\ell)=5$. 
\end{itemize}
This implies the desired claim (1). 

By using the similar argument, one can prove (2). 
The claim (3) is obvious from Lemma~\ref{lem_type_fiber}(3).  
$\Box$

\begin{lem}
\label{lem_dilatation_kl}
Suppose that $0 < \ell< k$. 
Then 
$$\lambda(k \mathfrak{a} \pm  \ell \mathfrak{b}) = \lambda(k \mathfrak{p} \pm  \ell \mathfrak{q}) = \lambda(k \mathfrak{r} \pm \ell \mathfrak{s}) = \lambda_{(k,\ell)}.$$
\end{lem}

\noindent
{\it Proof.} 
We use Theorem~\ref{thm_Fried-Oertel-Poly}. 
The dilatations of $ k \mathfrak{a} \pm  \ell \mathfrak{b}$ and $k \mathfrak{p} \pm \ell \mathfrak{q}$ 
are the largest real root of 
$$
f_{(2k \pm \ell, 2k \pm  2\ell, k \pm  2 \ell)}(t) = f_{(k, 2k \pm  2\ell,   \pm \ell)}(t) = 
(t^{k \pm  \ell}+1)(t^{2k} - t^{k+\ell} - t^k - t^{k-\ell}+1).
$$
The dilatation of $k \mathfrak{r} \pm \ell \mathfrak{s}$ is the largest real root of 
$$
f_{(k\pm\ell, k \mp \ell,  -k)}(t) 
=  (t^k+1)(t^{2k} - t^{k+\ell} - t^k - t^{k-\ell}+1).
$$
Since the absolute values of all roots of $t^{k \pm  \ell}+1$ and $t^k+1$ are equal to $1$,  
one finishes the proof. 
$\Box$
\medskip

By Proposition~\ref{prop_magic_ori} and (\ref{equation_3families}), we immediately obtain the following.

\begin{cor}
\label{cor_kl_orientable}
Suppose that $k >0$, $-k < \ell< k$ and $\gcd(k,\ell)=1$. 
\begin{description}
\item[(1)]  
The monodromy of the fibration associated to $ k \mathfrak{a}+ \ell \mathfrak{b}$ on $N$ is orientable if and only if $k$ is odd and $\ell$ is even. 
\item[(2)] 
The monodromy of the fibration associated to $k \mathfrak{p} + \ell \mathfrak{q}$ on $N$ is orientable if and only if $k$ is even and $\ell$ is odd. 
\item[(3)] 
The monodromy of the fibration associated to $k \mathfrak{r} + \ell \mathfrak{s} $ on $N$ is orientable if and only if both $k$ and  $\ell$ are odd. 
\end{description}
\end{cor}

The following can be obtained from Proposition~\ref{prop_sing-data}

\begin{cor}
\label{cor_sing_2cusps}
Suppose that $k >0$, $-k < \ell< k$ and  $\gcd(k,\ell)=1$. 
\begin{description}
\item[(1)]
The singularity data of the  monodromy  of the fibration associated to 
 $ k \mathfrak{a}+ \ell \mathfrak{b}$ on $N$ is given by 
$$(\underbrace{\tfrac{2k+ \ell}{\gcd (2k+\ell,5)} -2, \cdots, \tfrac{2k+ \ell}{\gcd (2k+\ell,5)} -2}_{\gcd (2k+\ell,5)},  
\underbrace{\tfrac{2k-\ell}{\gcd (5,k+2 \ell)} -2, \cdots, \tfrac{2k-\ell}{\gcd (5,k+2 \ell)}}_{\gcd (5,k+2 \ell)}).$$
\item[(2)] 
The singularity data of the  monodromy  of the fibration associated to 
$k \mathfrak{p} + \ell \mathfrak{q}$ on $N$  is given by 
$$(\underbrace{\tfrac{k}{\gcd (k,3)} -2, \cdots, \tfrac{k}{\gcd (k,3)} -2}_{\gcd (k,3)},  
\underbrace{\tfrac{3k}{\gcd (3,\ell)} -2, \cdots, \tfrac{3k}{\gcd (3,\ell)}-2}_{\gcd (3,\ell)}).$$
\item[(3)]
The singularity data of the monodromy  of the fibration associated to 
 $k \mathfrak{r} + \ell \mathfrak{s} $ on $N$  is given by 
$$(k+\ell -2 , k-\ell -2,   \underbrace{2, \cdots, 2}_{k}).$$
\end{description}
\end{cor}

\begin{rem}\ 
\label{rem_prong}
\begin{description}
\item[(1)] 
The stable foliation for the monodromy  of the fibration associated  to  $ k \mathfrak{a}+ \ell \mathfrak{b}$ (resp. $ k \mathfrak{p}+ \ell \mathfrak{q}$) 
has $2$ prongs on each boundary component on $T_{\beta}$. 
(Hence there is no singular leaf on $\partial_{\beta} F_{k \mathfrak{a}+ \ell \mathfrak{b}}$ (resp. $\partial_{\beta} F_{k \mathfrak{p}+ \ell \mathfrak{q}}$).) 
\item[(2)] 
The stable foliation for the monodromy  of the fibration associated to $ k \mathfrak{r}+ \ell \mathfrak{s}$ has $4$ prongs on each boundary component on $T_{\gamma}$. 
\end{description}
\end{rem}

\begin{lem}
\label{lem_1-prong_2cusps}
Suppose that $k >0$, $-k < \ell< k$ and $\gcd(k,\ell)=1$. 
\begin{description}
\item[(1)]  
The stable foliation for the  monodromy of the fibration associated to  $ k \mathfrak{a}+ \ell \mathfrak{b}$
 has a $1$ prong on a boundary component if and only if 
$(k,\ell) \in \{ (2, \pm 1), (3, \pm 1), (4, \pm 3)\}$. 
\item[(2)] 
The stable foliation for the  monodromy  of the fibration associated to $k \mathfrak{p} + \ell \mathfrak{q}$ has a $1$ prong on a boundary component  if and only if 
$(k,\ell) \in \{(1, 0), (3, \pm 1), (3, \pm 2)\}$. 
\item[(3)]
The stable foliation for the  monodromy of the fibration associated to $k \mathfrak{r} + \ell \mathfrak{s}$ has a $1$ prong on a boundary component  if and only if 
$k + \ell= 1$ or $k - \ell = 1$.  
\end{description}
\end{lem}

\noindent
{\it Proof.} 
(1) By Corollary~\ref{cor_sing_2cusps}, 
the stable foliation of the monodromy  for $ k \mathfrak{a}+ \ell \mathfrak{b}$ has a $1$ prong  if and only if 
 $\tfrac{2k+ \ell}{\gcd (2k+\ell,5)} = 1$ or $\tfrac{2k-\ell}{\gcd (5,k+2 \ell)}=1$. 
 
 Suppose that $\tfrac{2k+ \ell}{\gcd (2k+\ell,5)} = 1$. 
 Clearly $\gcd (2k+\ell,5)= 1$ or $5$.  
 If $\gcd (2k+\ell,5)=1$, then $2k+\ell=1$. 
 Since $-k < \ell < k$, one has $-k < -2k+1 < k$ which implies that $\tfrac{1}{3}< k< 1$. 
 This does not occur since $k$ is an integer. 
 If $\gcd (2k+\ell,5)=5$, then $2k+\ell = 5$. 
 Since $-k < \ell < k$, one has 
 $-k < 5-2k< k$ which implies $\tfrac{5}{3} < k < 5$. 
Hence $(k,\ell) \in \{(2,1), (3,-1), (4,-3)\}$. 

Suppose that $\tfrac{2k-\ell}{\gcd (5,k+2 \ell)}=1$. 
In this case one sees that 
$(k,\ell) \in \{(2,-1), (3,1), (4,3)\}$. 
This completes the proof of (1). 

By using the same argument one can prove (2),(3). 
$\Box$

\begin{rem}
\label{rem_not-used}
Suppose that $k >0$, $-k < \ell< k$ and $\gcd(k,\ell)=1$. 
By Lemma~\ref{lem_1-prong_2cusps} and Corollary~\ref{cor_magic-list}, we see that: 
\begin{description}
\item[(1)] 
$N(\tfrac{-3k-4\ell}{2k+\ell}, \tfrac{-3}{2}, \tfrac{-4k-3\ell}{k+2\ell})$ is non-hyperbolic 
for $(k,\ell) \in \{ (2, \pm 1), (3, \pm 1), (4, \pm 3)\}$. 
Otherwise it is hyperbolic. 
\item[(2)]
$N(\tfrac{-2k-3\ell}{k}, \tfrac{-1}{2}, \tfrac{-3k-2\ell}{\ell})$ is non-hyperbolic 
for $(k,\ell) \in \{(1, 0), (3, \pm 1), (3, \pm 2)\}$. 
Otherwise it is hyperbolic. 
\item[(3)] 
Suppose that  $k+\ell=1$ or $k-\ell=1$. 
Then we have the following. 
$N(\tfrac{\ell}{k+\ell}, \tfrac{-\ell}{k-\ell}, 2)$ is non-hyperbolic  if 
$(k,\ell) \in \{(2, \pm 1), (3, \pm 2), (4, \pm3)\}$ and 
it is  hyperbolic  if $(k,\ell) \notin \{(2, \pm 1), (3, \pm 2), (4, \pm3)\}$. 
Suppose that $k + \ell \ne  1$ and  $k - \ell \ne 1$. 
Then $N(\tfrac{\ell}{k+\ell}, \tfrac{-\ell}{k-\ell}, 2)$ is hyperbolic. 
\end{description}
\end{rem}

\begin{figure}
\begin{center}
\includegraphics[width=4in]{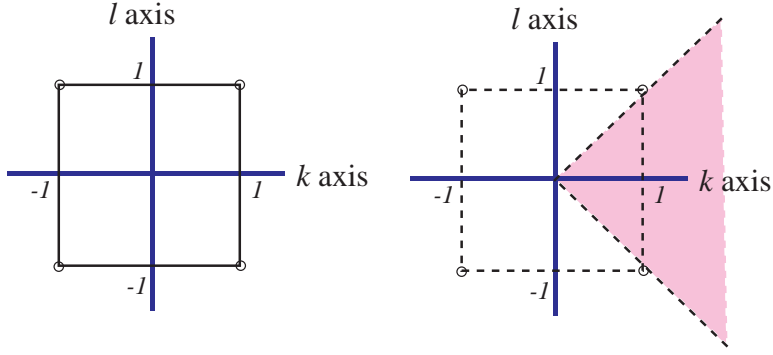}
\caption{(left) Thurston norm ball.  
(right) open cone $int(C_{\Delta(r)})$ [shaded region].}
\label{fig_ball}
\end{center}
\end{figure}

For each of  $N(\tfrac{-3}{2})$, $N(\tfrac{-1}{2})$ and $N(2)$, 
its Thurston norm ball of radius $2$ is a rectangle with vertices $(k, \ell)=(\pm 1, \pm1)$ 
illustrated in Figure~\ref{fig_ball}. 
By using  (\ref{equation_norm}) 
we  see that for $k , \ell \in {\Bbb R}$, 
\begin{equation}
\label{equation_FillingNorm}
\|k \overline{\mathfrak{a}}+ \ell \overline{\mathfrak{b}}\| = \|k \overline{\mathfrak{p}}+ \ell \overline{\mathfrak{q}}\| = 
\|k \overline{\mathfrak{r}}+ \ell \overline{\mathfrak{s}}\| = \max\{2|k|, 2 |\ell|\}.
\end{equation}

The following lemma asserts that fibered faces for $N(\tfrac{-3}{2})$ and $N(\tfrac{-1}{2})$ has a symmetry. 
Thus, for the study of monodromies on fibrations on  $N(\tfrac{-3}{2})$ (resp. $N(\tfrac{-1}{2})$), it is enough to consider the open cone over an arbitrary picked fibered face. 

\begin{lem}
\label{lem_symmetry-2cusp}
Suppose that $k >0$ and $- k < \ell < k$. 
\begin{description}
\item[(1)] 
$F_{ k \overline{\mathfrak{a}}+ \ell \overline{\mathfrak{b}}} = F_{ \ell \overline{\mathfrak{a}}+ k \overline{\mathfrak{b}}}$ and 
$\lambda(k \overline{\mathfrak{a}}+ \ell \overline{\mathfrak{b}}) = \lambda(\ell \overline{\mathfrak{a}}+ k \overline{\mathfrak{b}})$. 
\item[(2)] 
$F_{ k \overline{\mathfrak{p}}+ \ell \overline{\mathfrak{q}}} = F_{ \ell \overline{\mathfrak{p}}+ k \overline{\mathfrak{q}}}$ and 
$\lambda(k \overline{\mathfrak{p}}+ \ell \overline{\mathfrak{q}}) = \lambda(\ell \overline{\mathfrak{p}}+ k \overline{\mathfrak{q}})$. 
\end{description}
\end{lem}

\noindent
{\it Proof.} 
See Lemma~\ref{lem_symmetry} and Remark~\ref{rem_prong}(1). 
$\Box$
\medskip

Let us fix  an open cones 
\begin{eqnarray*}
int(C_{\Delta(-3/2)})&=& \{ k \overline{\mathfrak{a}}+ \ell \overline{\mathfrak{b}}\ |\ k > 0, \ -k < \ell < k\} 
\subset H_2(N(\tfrac{-3}{2}), \partial N(\tfrac{-3}{2}); {\Bbb R}), 
\\
int(C_{\Delta(-1/2)}) &=& \{k \overline{\mathfrak{p}} + \ell \overline{\mathfrak{q}}\ |\ k > 0, \ -k < \ell < k\} \subset H_2(N(\tfrac{-1}{2}), \partial N(\tfrac{-1}{2}); {\Bbb R}), 
\\
int(C_{\Delta(2)})&=& \{k \overline{\mathfrak{r}} + \ell \overline{\mathfrak{s}}\ |\ k > 0, \ -k < \ell < k\} \subset H_2(N(2), \partial N(2); {\Bbb R}).
\end{eqnarray*}

\noindent
Lemmas~\ref{lem_sym_fib} and \ref{lem_sym_dil} tell us that 
it is enough to consider the fibered classes of  $ int(C_{\Delta(r)})$ for  $0 < \ell< k$.

\begin{lem}
\label{lem_sym_fib}
Suppose that $0 < \ell< k$ and $\gcd(k,\ell)= 1$. 
\begin{description}
\item[(1)] 
$F_{ k \overline{\mathfrak{a}} \pm \ell \overline{\mathfrak{b}}}= \varSigma_{k-2,6}$  if $\gcd (2k+\ell,5)=5$ or $\gcd (5,k+2\ell)=5$. 
Otherwise $F_{ k \overline{\mathfrak{a}}+ \ell \overline{\mathfrak{b}}}= \varSigma_{k,2}$. 
\item[(2)] 
$F_{k \overline{\mathfrak{p}} \pm \ell \overline{\mathfrak{q}} }= \varSigma_{k-1, 4}   $ if $\gcd(k,3)=3$ or $\gcd(3,\ell)=3$. 
Otherwise $F_{k \overline{\mathfrak{p}} + \ell \overline{\mathfrak{q}}}= \varSigma_{k,2}  $. 
\item[(3)] 
$F_{k \overline{\mathfrak{r}} \pm \ell \overline{\mathfrak{s}} }=\varSigma_{k,2}$. 
\end{description}
\end{lem}

\noindent
{\it Proof.} 
The number of the components of $\partial_{\beta} F_{k \mathfrak{a} + \ell \mathfrak{b}}$ equals $k+ \ell$. 
By Lemma~\ref{lem_type_fiber}(1), we have the desired claim (1). 
One can prove (2),(3) by using Lemma~\ref{lem_type_fiber}(2),(3) respectively. 
$\Box$

\begin{lem}
\label{lem_sym_dil}
Suppose that $0 < \ell< k$. 
Then 
$$\lambda(k \overline{\mathfrak{a}} \pm \ell \overline{\mathfrak{b}})
= \lambda(k \overline{\mathfrak{p}} \pm \ell \overline{\mathfrak{q}})
=\lambda(k \overline{\mathfrak{r}} \pm \ell \overline{\mathfrak{s}})= \lambda_{(k,\ell)}.$$
\end{lem}

\noindent
{\it Proof.} 
See  Lemma~\ref{lem_dilatation_kl} and Remark~\ref{rem_prong}. 
$\Box$

\begin{prop}
\label{prop_mini-ray_2cusps}
Let $r \in \{ \tfrac{-3}{2}, \tfrac{-1}{2}, 2\}$. 
The minimum of  $\mathrm{Ent}: int(C_{\Delta(r)}) \rightarrow {\Bbb R}$  equals $2 \log \lambda_{(1,0)}=2 \log \bigl(\tfrac{3+ \sqrt{5}}{2} \bigr) $. 
The minimizer is given by  
$\overline{ \mathfrak{a}}$ if $r=  \tfrac{-3}{2}$, 
$\overline{ \mathfrak{p}}$ if $r=  \tfrac{-1}{2}$ and 
$ \overline{ \mathfrak{r}}$ if  $r=  2$.
\end{prop}

\noindent
{\it Proof.}
Recall that  $\mathrm{Ent}$ is constant on each ray thorough the origin 
and it has the minimum at a unique ray.  
By Lemma~\ref{lem_sym_dil} and (\ref{equation_FillingNorm}), the ray which reaches the minimum must satisfy $\ell=0$. 
As a representative which lies on  such a ray, we take $k=1$ and $\ell=0$. 
This together with Theorem~\ref{thm_Fried-Oertel-Poly} implies that 
$\min \mathrm{Ent}= 2 \log \lambda_{(1,0)}=2 \log \bigl(\tfrac{3+ \sqrt{5}}{2} \bigr)$. 
$\Box$ 

\begin{rem}
The monodromies of the fibrations associated to 
$\overline{ \mathfrak{a}} $ (resp. $\overline{ \mathfrak{p}}$) on $N(\tfrac{-3}{2})$ (resp. $N(\tfrac{-1}{2})$)  are intriguing examples. 
\begin{description}
\item[(1)] 
$N(\tfrac{-3}{2})$ admits a fiber of genus $1$ with $2$ boundary components 
corresponding to $\overline{ \mathfrak{a}} $. 
Its monodromy $\Phi: \varSigma_{1,2} \rightarrow \varSigma_{1,2}$ has $2$ prongs on each boundary component. 
Thus $\Phi$ extends to the monodromy $\overline{\Phi}: \varSigma_{1,1} \rightarrow \varSigma_{1,1}$ 
of the fibration on $N(\tfrac{-3}{2}, \tfrac{-3}{2})$ (which is the figure-$8$ knot sister manifold, see \cite[Table~A.2]{MP}) 
by capping the  boundary component  of $\varSigma_{1,2}$ on $T_{\alpha}$. 
It is well-known that $\overline{\Phi}$ realizes the minimal dilatation $ \tfrac{3+ \sqrt{5}}{2} $ among pseudo-Anosovs 
on $\varSigma_{1,1}$. 
\item[(2)] 
$N(\tfrac{-1}{2})$ admits a fiber of genus $0$ with $4$ boundary components 
corresponding to $\overline{ \mathfrak{p}}$. 
Its monodromy $\Phi$ fixes a boundary component, and hence it can be considered that $\Phi$ is a pseudo-Anosov homeomorphism on a $3$-punctured disk $D_3$. 
This monodromy $\Phi$ realizes the minimal dilatation $\delta(D_3)= \tfrac{3+ \sqrt{5}}{2} $ 
among pseudo-Anosovs on $D_3$. 
\end{description}
\end{rem}

\subsection{Property of algebraic integers $\lambda_{(k,\ell)}$}

\begin{lem}
\label{lem_easy-monotonicity}
Suppose that $1 < \ell+1 < k$ and $\gcd(k,\ell)=1$. 
Then 
$\lambda_{(k+1,\ell)}< \lambda_{(k,\ell)}< \lambda_{(k,\ell+1)}$. 
\end{lem}

\noindent
{\it Proof.} 
The ray which attains the minimum of $ \mathrm{Ent}: int(C_{\Delta(r)}) \rightarrow {\Bbb R}$ 
satisfies $\ell=0$. 
Recall that $\tfrac{1}{\mathrm{ent}(\cdot)}: int(C_{\Delta(r)}(\Bbb Q)) \rightarrow {\Bbb R}$ is strictly concave. 
Thus one has 
$$\log \lambda_{(k, \tfrac{k \ell}{k+1})}< \log \lambda_{(k,\ell)} < \log \lambda_{(k, \ell+1)}.$$
The inequality $\log \lambda_{(k+1,\ell)}< \log \lambda_{(k,\ell)}$ holds 
since 
$$\log \lambda_{(k+1,\ell)}  =\mathrm{ent}((k+1)\overline{\mathfrak{a}}+ \ell \overline{\mathfrak{b}}) 
= \tfrac{k}{k+1} \mathrm{ent}(k \overline{\mathfrak{a}}+    \tfrac{k \ell}{k+1}  \overline{\mathfrak{b}}) = \tfrac{k}{k+1} \log \lambda_{(k, \tfrac{k \ell}{k+1})} 
<  \log \lambda_{(k, \tfrac{k \ell}{k+1})} .$$
Hence $\log \lambda_{(k+1,\ell)}< \log \lambda_{(k,\ell)}< \log \lambda_{(k,\ell+1)}$. 
$\Box$ 

\begin{lem}
\label{lem_asymp_roots}
For any fixed $\ell>0$, 
$$\lim_{k \to \infty}k  \log \lambda_{(k,\ell)} =\log \lambda_{(1,0)}= \log (\tfrac{3 + \sqrt{5}}{2}).$$
\end{lem}

\noindent
{\it Proof.} 
The ray through $k \overline{\mathfrak{a}} + \ell \overline{\mathfrak{b}}$ from the origin goes to the ray through $ \overline{\mathfrak{a}} $ 
if $k$ goes to $\infty$. 
Hence 
$$\lim_{k \to \infty} \mathrm{Ent}(k \overline{\mathfrak{a}} + \ell \overline{\mathfrak{b}}) 
= \lim_{k \to \infty} 2k \log \lambda_{(k,\ell)}
= \mathrm{Ent}(\overline{\mathfrak{a}}) =2 \log \lambda_{(1,0)}.$$ 
This completes the proof. 
$\Box$ 

\begin{prop}
\label{prop_monotone}
If $\lambda_{(k+1, \ell)} < \lambda_{(k,1)}$ for some $k \ge \ell \ge 2$, then 
$ \lambda_{(k+2, \ell)} < \lambda_{(k+1,1)}$. 
\end{prop}

\begin{figure}
\begin{center}
\includegraphics[width=3in]{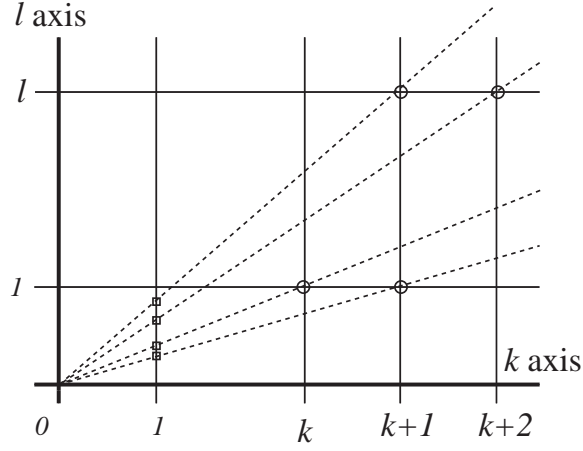}
\caption{four boxes $\Box$  (from the bottom to the top) on the line $k=1$ indicate $p_1,p_2,p_3$ and $p_4$.}
\label{fig_4point}
\end{center}
\end{figure}

\noindent
{\it Proof.} 
We denote the homology class $k \overline{\mathfrak{a}}+ \ell \overline{\mathfrak{b}}$ by $(k,\ell)$. 
Since  $ \mathrm{Ent}$ is constant on each ray thorough the origin, 
$k \mathrm{ent}(k,\ell)=  \mathrm{ent}(1, \tfrac{\ell}{k})$. 
One takes $4$ points 
$$p_1= (1, \tfrac{1}{k+1}),\ 
p_2= (1, \tfrac{1}{k}),\ 
p_3= (1, \tfrac{\ell}{k+2}),\ 
p_4= (1, \tfrac{\ell}{k+1}),$$ 
see Figure~\ref{fig_4point}. 
We have $\tfrac{1}{k+1} <   \tfrac{1}{k} < \tfrac{2}{k+2} \le  \tfrac{\ell}{k+2} <  \tfrac{\ell}{k+1}$. 
Let us set $t,t'$ and $c$ as follows. 
\begin{equation}
\label{equation_dainyu}
\begin{split}
0 < t&=\frac{|p_3-p_2|}{|p_4-p_2|}= \frac{(k+1)(k\ell-k-2)}{(k+2)(k \ell-k-1)}<1, 
\\
0 < t&'=\frac{|p_3-p_2|}{|p_3-p_1|}=  \frac{|p_4-p_2|}{|p_3-p_1|}t   = \frac{(k+2)(k \ell-k-1)}{k(k \ell-k+\ell-2)} t<1, 
\\
1< c&=  \frac{(k+2)(k \ell-k-1)}{k(k \ell-k+\ell-2)}. 
\end{split}
\end{equation}
(Hence $t'= ct$.)
Then 
\begin{eqnarray*}
|p_3-p_2|:| p_4-p_3|&=& t: 1-t, 
\\
|p_3-p_2|: |p_2-p_1|&=& ct: 1-ct.
\end{eqnarray*}
These ratios together with Theorem~\ref{thm_strictly-concave} imply that 
\begin{eqnarray}
\label{equation_concave1}
\frac{1}{(k+2) \mathrm{ent}(k+2,\ell)} > (1-t) \frac{1}{k \mathrm{ent}(k,1)}+ t \frac{1}{(k+1)\mathrm{ent}(k+1, \ell)},
\\
\label{equation_concave2}
\frac{1}{k \mathrm{ent}(k,1)} > ct \frac{1}{(k+1) \mathrm{ent}(k+1,1)}+ (1-ct) \frac{1}{(k+2) \mathrm{ent}(k+2, \ell)}.
\end{eqnarray}
By (\ref{equation_concave1}) and by the assumption $\mathrm{ent}(k,1) > \mathrm{ent}(k+1, \ell)$, 
\begin{eqnarray*}
\frac{1}{(k+2) \mathrm{ent}(k+2,\ell)} 
&>& (1-t) \frac{1}{k \mathrm{ent}(k,1)}+ t \frac{1}{(k+1)\mathrm{ent}(k+1, \ell)}
\\
&>&  (1-t) \frac{1}{k \mathrm{ent}(k,1)}+ t \frac{1}{(k+1)\mathrm{ent}(k, 1)}
\\
&=& \frac{k+1-t}{k+1} \frac{1}{k \mathrm{ent}(k,1)} 
\\
&>& \frac{k+1-t}{k+1} \Big\{ ct \frac{1}{(k+1) \mathrm{ent}(k+1,1)}+ (1-ct) \frac{1}{(k+2) \mathrm{ent}(k+2, \ell)}\Big\}. 
\end{eqnarray*}
The last inequality is given by  (\ref{equation_concave2}). 
Hence 
\begin{eqnarray*}
\Big\{ \frac{1}{k+2}\ - \frac{(k+1-t)(1-ct)}{(k+1)(k+2)}\Big\} \frac{1}{\mathrm{ent}(k+2, \ell)} > 
\frac{(k+1-t)ct}{(k+1)^2} \frac{1}{\mathrm{ent}(k+1,1)}, 
\end{eqnarray*}
which gives, by calculation, 
\begin{eqnarray*}
\frac{(k+1-t)c+1}{k+2}  \frac{1}{ \mathrm{ent}(k+2,\ell)} 
> 
\frac{(k+1-t)c}{k+1} \frac{1}{\mathrm{ent}(k+1,1)}.
\end{eqnarray*}
Thus 
\begin{eqnarray*}
\mathrm{ent}(k+2,\ell) < \Big\{\frac{k+1}{(k+1-t)c}\Big\} \Big\{ \frac{(k+1-t)c+1}{k+2} \Big\}\mathrm{ent}(k+1,1). 
\end{eqnarray*}
For the proof of the claim 
it is enough to verify the equality 
$ \Big\{\frac{k+1}{(k+1-t)c}\Big\} \Big\{ \frac{(k+1-t)c+1}{k+2} \Big\} =1$. 
Clearly, 
\begin{eqnarray*}
&\ & 
\Big\{\frac{k+1}{(k+1-t)c}\Big\} \Big\{ \frac{(k+1-t)c+1}{k+2} \Big\} =1 
\\
&\Leftrightarrow& 
(k+1) \{ (k+1-t)c+1\} = (k+2)(k+1-t)c
\\
&\Leftrightarrow&  
(k+1)(k+1-t)c + k+1 = (k+2)(k+1-t)c
\\
&\Leftrightarrow&  
k+1 = (k+1-t)c 
\\
&\Leftrightarrow&   
c= \frac{k+1}{k+1-t}. 
\end{eqnarray*}
One can verify the last equality $c= \frac{k+1}{k+1-t}$ by substituting the constants $t$ and $c$ given by (\ref{equation_dainyu}). 
$\Box$

\subsection{Fibrations of manifolds obrained from $N(\tfrac{-1}{2}), N(\tfrac{-3}{2})$ and $N(2)$ by Dehn filling  two cusps}

As a consequence of  Lemma~\ref{lem_1-prong_2cusps}, we see the following.  

\begin{rem}
\label{rem_extended-pA}
If $(k,\ell) \notin \{ (2, \pm 1), (3, \pm 1), (4, \pm 3)\}$, then 
the monodromy $ \Phi_{k \mathfrak{a}+ \ell \mathfrak{b}}:  F_{k \mathfrak{a}+ \ell \mathfrak{b}} \rightarrow F_{k \mathfrak{a}+ \ell \mathfrak{b}}$ 
of the fibration  associated to $k \mathfrak{a}+ \ell \mathfrak{b}$ on $N$ 
extends to the  monodromy 
$ \overline{\Phi}_{k \mathfrak{a}+ \ell \mathfrak{b}}: \overline{F}_{k \mathfrak{a}+ \ell \mathfrak{b}} \rightarrow  \overline{F}_{k \mathfrak{a}+ \ell \mathfrak{b}} $ 
of the fibration on $N(\tfrac{-3k-4\ell}{2k+\ell}, \tfrac{-3}{2}, \tfrac{-4k-3\ell}{k+2\ell})$ with the dilatation $\lambda_{(k,\ell)} (= \lambda(\Phi_{k \mathfrak{a}+ \ell \mathfrak{b}}))$. 
Similarly, if $(k,\ell) \notin \{(1, 0), (3, \pm 1), (3, \pm 2)\}$ (resp. if $k + \ell \ne 1$ and  $k - \ell \ne 1$), 
then the monodromy $\Phi_{k \mathfrak{p}+ \ell \mathfrak{q}}$ (resp. $\Phi_{k \mathfrak{r}+ \ell \mathfrak{s}}$) of the fibration 
associated to  $k \mathfrak{p}+ \ell \mathfrak{q}$ (resp. $k \mathfrak{r}+ \ell \mathfrak{s}$) on $N$  
extends to the  monodromy $\overline{\Phi}_{k \mathfrak{p}+ \ell \mathfrak{q}}$ 
(resp.  $\overline{\Phi}_{k \mathfrak{r}+ \ell \mathfrak{s}}$) 
of the fibration on  $N(\tfrac{-2k-3\ell}{k}, \tfrac{-1}{2}, \tfrac{-3k-2\ell}{\ell})$ 
(resp. $N(\tfrac{\ell}{k+\ell}, \tfrac{-\ell}{k-\ell}, 2)$) with the dilatation $\lambda_{(k,\ell)}$. 
\end{rem}

\noindent
Let $\overline{\phi}_{k \mathfrak{a}+ \ell \mathfrak{b}}$, $\overline{\phi}_{k \mathfrak{p}+ \ell \mathfrak{q}}$ and $\overline{\phi}_{k \mathfrak{r}+ \ell \mathfrak{s}}$ 
be elements of $\mathrm{Mod}(\varSigma_g)$ containing  
$\overline{\Phi}_{k \mathfrak{a}+ \ell \mathfrak{b}}$, $\overline{\Phi}_{k \mathfrak{p}+ \ell \mathfrak{q}}$ and $\overline{\Phi}_{k \mathfrak{r}+ \ell \mathfrak{s}}$ 
as a representative. 

\begin{prop}
\label{prop_volume-conv}
For any fixed integer $\ell>0$, we have the following. 
\begin{description}
\item[(1)]

$\displaystyle \lim_{\substack{k  \to \infty \\ \gcd(k,\ell)=1}} \mathrm{vol}(\overline{\phi}_{k \mathfrak{a} + \ell \mathfrak{b}}) = \mathrm{vol}(N(\tfrac{-3}{2})) \approx 3.66386$. 
\item[(2)] 
$\displaystyle \lim_{\substack{k \to  \infty \\ \gcd(k,\ell)=1}} \mathrm{vol}(\overline{\phi}_{k \mathfrak{p} + \ell \mathfrak{q}}) = \mathrm{vol}(N(\tfrac{-1}{2}))  \approx 4.05977$. 
\item[(3)] 
$\displaystyle \lim_{\substack{k  \to  \infty\\  \gcd(k,\ell)=1}} \mathrm{vol}(\overline{\phi}_{k \mathfrak{r} + \ell \mathfrak{s}}) = \mathrm{vol}(N(2))  \approx  4.41533$. 
\end{description}
\end{prop}

\noindent
{\it Proof.} 
We will prove the claim (1). 
The proof of claims (2),(3) is similar. 
The mapping torus ${\Bbb T}(\overline{\phi}_{k \mathfrak{a} + \ell \mathfrak{b}}) $ is homeomorphic to 
$N(\tfrac{-3k-4\ell}{2k+\ell}, \tfrac{-3}{2}, \tfrac{-4k-3\ell}{k+2\ell})$. 
Since $\gcd(-3k-4\ell, 2k+\ell)$ (resp. $\gcd(-4k-3\ell, k+2\ell)$) is either $1$ or $5$, 
the two points 
$$(\tfrac{-3k-4\ell}{\gcd(-3k-4\ell, 2k+\ell)}, \tfrac{2k+\ell}{\gcd(-3k-4\ell, 2k+\ell)})\in {\Bbb R}^2, 
(\tfrac{-4k-3\ell}{\gcd(-4k-3\ell, k+2\ell)}, \tfrac{k+2\ell}{\gcd(-4k-3\ell, k+2\ell)}) \in {\Bbb R}^2$$
 tend to $\infty$ as $k$ tends to $\infty$. 
Thurston's hyperbolic Dehn surgery theorem (see \cite{Thurston}) implies 
the volume of $N(\tfrac{-3k-4\ell}{2k+\ell}, \tfrac{-3}{2}, \tfrac{-4k-3\ell}{k+2\ell})$ converges to $\mathrm{vol}(N(\tfrac{-3}{2}))$ as $k$ tends to $\infty$. 
$\Box$
\medskip

\noindent
{\it Proof of Theorem~\ref{thm_three}.} 
(Case $r = \tfrac{-3}{2}$.) 
For the proof of (1), first of all we  find a pair $(k(g),  \ell(g))= (g+ \tilde{k}(g), \ell(g))$ for each $g \ge 3$ satisfying the following: 
the both $ \tilde{k}(g) >0$, $\ell(g)>0$ are bounded, and 
the genus of  $F_{k(g)  \mathfrak{a}+  \ell(g)\mathfrak{b}}$ equals $g$. 
Next we  check that the stable foliation of  $\Phi_{k(g)  \mathfrak{a}+  \ell(g)\mathfrak{b}}$ has no $1$ prong 
on each boundary component of $F_{k(g)  \mathfrak{a}+  \ell(g)\mathfrak{b}}$. 
Then one can extend $\Phi_{k(g)  \mathfrak{a}+  \ell(g)\mathfrak{b}}$ to the pseudo-Anosov homeomorphism 
$\overline{\Phi}_{k(g)  \mathfrak{a}+  \ell(g)\mathfrak{b}}$ 
on a closed surface of genus $g$. 
This finishes the proof of (1). 
In fact, by Lemma~\ref{lem_asymp_roots} 
$$\lim_{g \to \infty} k(g) \log \lambda(\overline{\Phi}_{k(g)  \mathfrak{a}+  \ell(g)\mathfrak{b}})
= \lim_{g \to \infty} k(g) \log \lambda_{(k(g), \ell(g))} = \log (\tfrac{3+ \sqrt{5}}{2}).$$ 
On the other hand 
$$\lim_{g \to \infty} \log \lambda(\overline{\Phi}_{k(g)  \mathfrak{a}+  \ell(g)\mathfrak{b}})
= \lim_{g \to \infty} \tfrac{1}{g+  \tilde{k}(g)} \log ( \tfrac{3+ \sqrt{5}}{2}) = 0.$$
Thus one obtains 
\begin{eqnarray*}
\log( \tfrac{3+ \sqrt{5}}{2}) = \lim_{g \to \infty} k(g) \log \lambda(\overline{\Phi}_{k(g)  \mathfrak{a}+  \ell(g)\mathfrak{b}})
&=& \lim_{g \to \infty} (g+ \tilde{k}(g))  \log \lambda(\overline{\Phi}_{k(g)  \mathfrak{a}+  \ell(g)\mathfrak{b}})
\\
&=& \lim_{g \to \infty} g  \log \lambda(\overline{\Phi}_{k(g)  \mathfrak{a}+  \ell(g)\mathfrak{b}}) 
+ \lim_{g \to \infty}  \tilde{k}(g)  \log \lambda(\overline{\Phi}_{k(g)  \mathfrak{a}+  \ell(g)\mathfrak{b}})
\\
&=&  \lim_{g \to \infty} g   \log \lambda(\overline{\Phi}_{k(g)  \mathfrak{a}+  \ell(g)\mathfrak{b}})+0, 
\end{eqnarray*}
which implies (1). 

One sees that the genera of $F_{3  \mathfrak{a}+ 2 \mathfrak{b}}$ and $F_{4  \mathfrak{a}+  \mathfrak{b}}$ equal $3$ and $4$ respectively. 
If  $g \not\equiv 0 \pmod 5$ and $g  \ge 6$, the genus of $F_{g  \mathfrak{a}+ 5 \mathfrak{b}}$ equals $g$. 
In the case $g \equiv 2 \pmod 5$ and $g \ge 7$, the genus of $F_{g  \mathfrak{a}+  \mathfrak{b}}$ equals $g -2 \equiv 0 \pmod 5$. 
By Lemma~\ref{lem_1-prong_2cusps}, one has the desired equality (1).  

The claim (2) on the volume holds by Proposition~\ref{prop_volume-conv}(1). 
\medskip
\\
(Case $r = \tfrac{-1}{2}$.)
If $g \equiv 0, 1 \pmod 3$ and $g \ge 3$, 
the genus of $F_{(g+1)  \mathfrak{p}+ 3 \mathfrak{q}}$ equals $g$. 
If $g \equiv 2 \pmod 3$ and $g \ge 3$, 
the genus of  $F_{(g+1)  \mathfrak{p}+  \mathfrak{q}}$ equals $g$. 
By Lemma~\ref{lem_1-prong_2cusps} and Proposition~\ref{prop_volume-conv}(2), 
one obtains the claims (1),(2). 
\medskip
\\
(Case $r=2$.) 
The genus of $F_{g \mathfrak{r} +  \mathfrak{s}}$ equals $g$. 
By Lemma~\ref{lem_1-prong_2cusps} and Proposition~\ref{prop_volume-conv}(3), 
one obtains the claims (1),(2). 
$\Box$
\medskip

\begin{rem}
\label{rem_LT}
For $g \ge 4$ even,  there exists a $\varSigma_g$-bundle over the circle with the dilatation $\lambda_{(g,1)}$, 
which is obtained from the  extension of the monodromy of the fibration associated to $g \mathfrak{r} +  \mathfrak{s}$ on $N$. 
However the invariant foliation associated to this $\varSigma_g$-bundle over the circle is non-orientable, 
see Corollary~\ref{cor_kl_orientable} and Question~\ref{ques_LT}. 
(In  the case $g=2$, 
the monodromy of the fibration associated to $2 \mathfrak{r} +  \mathfrak{s}$ cannot extend to the {\it pseudo-Anosov} monodromy on a closed fiber, see Remark~\ref{rem_not-used}(3).) 
%The monodromy of the fibration associated to $g \mathfrak{r} +  \mathfrak{s}$ for $g$ even is not orientable, see Corollary~\ref{cor_kl_orientable}. 
%This tells us that for each $g \ge 4$ even, there exists a $\varSigma_g$-bundle over the circle  with non-orientable invariant foliation 
%whose dilatation equals $\lambda_{(g,1)}$ (cf. Question~\ref{ques_LT}), 
%which is obtained from $N(2)$ by Dehn filling two cusps. 
%(In  the case $g=2$, 
%the monodromy of the fibration associated to $2 \mathfrak{r} +  \mathfrak{s}$ cannot extend to the {\it pseudo-Anosov} monodromy on a closed fiber, see Remark~\ref{rem_not-used}(3).) 
\end{rem}

For  $r \in {\Bbb Q} $, 
 $\Lambda_g(r)$ (resp. $\Lambda_g^+(r)$) 
 is defined to be the set of dilatations of all $\varSigma_g$-bundles 
 (resp. all $\varSigma_g$-bundles with orientable invariant foliations) 
 which are obtained from $N(r)$ by Dehn filling two cusps along the boundary slopes of the fibers of $N(r)$. 
Recall that $\mathcal{U}$ and $\mathcal{U}^+$ are finite sets of fibered hyperbolic $3$-manifolds defined in the introduction. 

\begin{lem}
\label{lem_list-KT1} 
$N(2) \in \mathcal{U}^+$. 
\end{lem}

\noindent
{\it Proof.} 
One sees that the pseudo-Anosov 
$\overline{\phi}_{3 \mathfrak{r} +  \mathfrak{s}} \in \mathrm{Mod}(\varSigma_3)$ is orientable and it has the dilatation 
$\lambda_{(3,1)} (= \delta_3^+)$. 
Hence $\delta_3^+ \in \Lambda_3^+(2)$. 
$\Box$
\medskip

In the rest of this section, we consider the sets $ \Lambda_g^{(+)}(\tfrac{-1}{2})$ and $ \Lambda_g^{(+)}(\tfrac{-3}{2})$ mainly. 
We first recall the number  $\min \Lambda_g^{(+)}(\tfrac{-1}{2}) $.

\begin{prop}[\cite{Hironaka}] 
\label{prop_one-half}
Let $g \ge 3$. 
%Among the dilatations of monodromies on a closed fiber of genus $g$ obtained from Dehn fillings of $N(\tfrac{-1}{2})$, 
%the following describes the minimum. 
%Among homology classes $a$ of $H_2(N(\tfrac{-1}{2}), \partial N(\tfrac{-1}{2}); {\Bbb R})$ which induce monodromies $\overline{\Phi(a)}$ on a closed fiber of genus $g$ 
%by Dehn fillings of $N(\tfrac{-1}{2})$, the following describes a homology class which attains the minimum among dilatations of monodromies $\overline{\Phi(a)}$. 
\begin{description}
\item[(1)] 
$\lambda_{(g+1,3)} = \min \Lambda_g(\tfrac{-1}{2})$  if $g \equiv 0,1,3,4 \pmod 6$. 
\item[(2)] 
$\lambda_{(g+1,1)} = \min \Lambda_g(\tfrac{-1}{2})$ if $g \equiv 2,5 \pmod 6$. 
\end{description}
\end{prop}

\begin{prop}[\cite{Hironaka}] 
\label{prop_one-half-ori}
Let $g \ge 3$. 
\begin{description}
\item[(1)] 
$\lambda_{(g+1,3)} = \min \Lambda_g^+(\tfrac{-1}{2})$
if $g \equiv 1,3 \pmod 6$. 
\item[(2)] 
$\lambda_{(g,1)} = \min \Lambda_g^+(\tfrac{-1}{2})$ 
if $g \equiv 2,4 \pmod 6$. 
\item[(3)] 
$\lambda_{(g+1,1)} = \min \Lambda_g^+(\tfrac{-1}{2})$
if $g \equiv 5 \pmod 6$. 
\end{description}
\end{prop}

\begin{lem}[\cite{Hironaka}] 
\label{lem_list-H}
$N(\tfrac{-1}{2}) \in \mathcal{U} \cap \mathcal{U}^+$. 
\end{lem}

\noindent
{\it Proof.} 
One sees that 
$\overline{\phi}_{2 \mathfrak{p} +  \mathfrak{q}} \in \mathrm{Mod}(\varSigma_2)$ is an orientable pseudo-Anosov mapping class having dilatation 
$\lambda_{(2,1)} (=\delta_2 = \delta_2^+) $. 
Hence $\delta_2 = \delta_2^+  \in \Lambda_2(\tfrac{-1}{2}) \cap \Lambda^+_2(\tfrac{-1}{2}) $. 
$\Box$
\medskip

We turn to  $N(\tfrac{-3}{2})$. 
By Lemma~\ref{lem_type_fiber}(1), 
if $\lambda \in \Lambda_g(\tfrac{-3}{2})$, then 
$\lambda = \lambda_{(g+2, \ell)}$ for some $1 \le \ell < g+2$ or 
$\lambda= \lambda_{(g, \ell')}$ for some $1 \le \ell' < g$. 

It is easy to verify the following by a direct computation.

\begin{lem}
\label{lem_genus_k-2}
For integers $k$ and $\ell$, 
$\gcd (2k+\ell,5)=5$ or $\gcd (5,k+2\ell)=5$ 
if and only if $k$ and $\ell$ are either 
\begin{description}
\item[(1)] $ \ell \equiv 0 \pmod{5}$ and $k \equiv 0 \pmod{5}$,
\item[(2)] $\ell \equiv 1 \pmod{5}$ and $k \equiv 2,3 \pmod{5}$, 
\item[(3)] $\ell \equiv 2 \pmod{5}$ and $k \equiv 1,4 \pmod{5}$, 
\item[(4)] $\ell \equiv 3 \pmod{5}$ and $k \equiv 1,4 \pmod{5}$, or 
\item[(5)] $\ell \equiv 4 \pmod{5}$ and $k \equiv 2,3 \pmod{5}$. 
\end{description}
\end{lem}

\noindent
We compute $\min \Lambda_g(\tfrac{-3}{2})$ for $g  \equiv 0,1,3, 5,6, 7,8,9 \pmod {10}$.

\begin{prop}
\label{prop_bound_KT1}
\begin{description}
\item[(1)] 
$\lambda_{(g+2,1)}  =\min \Lambda_g(\tfrac{-3}{2}) < \min \Lambda_g(\tfrac{-1}{2}) $  if $g \equiv 0,1,5,6 \pmod {10}$ and $g \ge 5$. 
\item[(2)] 
$\lambda_{(g+2,2)} =\min \Lambda_g(\tfrac{-3}{2}) < \min \Lambda_g(\tfrac{-1}{2}) $ if $g \equiv 7,9 \pmod {10}$ and $g \ge 7$. 
\item[(3)] 
$\lambda_{(g,2)} =\min \Lambda_g(\tfrac{-3}{2}) > \min \Lambda_g(\tfrac{-1}{2}) $ if $g \equiv 3 \pmod {10} $ and $g \ge 3$. 
\item[(4)] 
Let $g \equiv 8 \pmod{10}$ and $g \ge 8$. 
\begin{description}
\item[(i)] 
 $\lambda_{(g,3)} =\min \Lambda_g(\tfrac{-3}{2}) > \min \Lambda_g(\tfrac{-1}{2}) $ if $g \equiv 8,28 \pmod{30}$, 
 \item[(ii)] 
 $\lambda_{(g,5)} =\min \Lambda_g(\tfrac{-3}{2}) > \min \Lambda_g(\tfrac{-1}{2}) $  if $g \equiv 18 \pmod{30}$. 
\end{description}
\end{description}
\end{prop}

\noindent
{\it Proof}.
(1) 
If $k \equiv 2,3 \pmod{5}$, then 
$\gcd(2k+1,5)= 5$ or $\gcd(5,k+2)=5$. 
We set $k= g+2$. (Hence $g  \equiv 0,1 \pmod{5}$ or equivalently $g \equiv 0,1,5,6 \pmod {10}$.)  
The genus of $F_{(g+2) \mathfrak{a}+ \mathfrak{b}}$ is equal to $g$ by Lemma~\ref{lem_type_fiber}(1), and hence 
$\lambda_{(g+2,1)} \in  \Lambda_g(\tfrac{-3}{2})$ by Remark~\ref{rem_extended-pA}. 
One can check that $\lambda_{(g+2,1)}$ attains $\min \Lambda_g(\tfrac{-3}{2})$ by Lemma~\ref{lem_easy-monotonicity}. 
In fact for any $g >1$, $1 \le \ell' < g$ and $1 \le \ell < g+2$, it follows that 
$$\lambda_{(g+2,1)} < \lambda_{(g+1,1)}< \lambda_{(g,1)} \le \lambda_{(g, \ell')} \ \mbox{and}\ 
\lambda_{(g+2,1)} < \lambda_{(g+2,\ell)}.$$
Thus $\lambda_{(g,1)} = \min \Lambda_g(\tfrac{-3}{2})$. 
%Thus $\lambda_{(g+2,1)} \le \lambda$ for any $\lambda \in  \Lambda_g(\tfrac{-3}{2})$. 

By Proposition~\ref{prop_one-half}, 
the lower and upper bound of $\min \Lambda_g(\tfrac{-1}{2})$ is given by 
\begin{equation}
\label{equation_bound-one-half}
\lambda_{(g+1,1)} \le \min \Lambda_g(\tfrac{-1}{2}) \le \lambda_{(g+1,3)} \hspace{2mm}\mbox{for\ any\ }g.
\end{equation}
Since $\lambda_{(g+2,1)}< \lambda_{(g+1,1)}$, one obtains the inequality 
$\min \Lambda_g(\tfrac{-3}{2}) < \min \Lambda_g(\tfrac{-1}{2}) $. 
\medskip
\\
(2) 
If $k \equiv 1,4 \pmod{5}$, then 
$\gcd(2k+2,5) = 5$ or $\gcd(5, k+4)= 5$. 
We set $k= g+2$. (Hence $g \equiv 2,4 \pmod{5}$.) 
Suppose that $\gcd(g+2,2)=1$. 
Then $\lambda_{(g+2,2)} \in \Lambda_g(\tfrac{-3}{2})$ and 
$g  \equiv 7,9 \pmod{10}$ since $g$ must be odd. 
One sees that $\lambda_{(g+2,1)} \notin \Lambda_g(\tfrac{-3}{2})$ 
since $\gcd(2k+1,5)= 1$ and  $\gcd(5,k+2)=1$. 
For any $g >1$ and $1 \le \ell < g$, it follows that  
$\lambda_{(g+1,1)}< \lambda_{(g,1)} \le \lambda_{(g, \ell)}$. 
On the other hand 
$$\lambda_{(5,2)} \approx 1.23039 < \lambda_{(4,1)} \approx 1.28064$$ 
%$$\lambda_{(3,2)} \approx  1.50614< \lambda_{(2,1)}\approx 1.72208$$ 
and by Proposition~\ref{prop_monotone}, 
one has $\lambda_{(g+2, 2)} < \lambda_{(g+1,1)}$ holds for any $g \ge 3$. 
Thus $\lambda_{(g+2,2)} $ attains $\min  \Lambda_g(\tfrac{-3}{2})$. 
The inequality $\min \Lambda_g(\tfrac{-3}{2}) < \min \Lambda_g(\tfrac{-1}{2}) $ holds by (\ref{equation_bound-one-half}). 
\medskip
\\
(3),(4) 
Suppose that $g \equiv 3 \pmod{5}$, that is $g \equiv 3,8 \pmod{10}$. 
One observes that the genus of $F_{(g+2) \mathfrak{a}+ \ell \mathfrak{b}}$ equals $g+2$ whenever $\gcd(g+2, \ell)=1$. 
Hence if $\lambda \in  \Lambda_g(\tfrac{-3}{2})$, then $\lambda= \lambda_{(g, \ell)}$ for some $1 \le \ell < g$. 
Suppose that $g \equiv 3 \pmod {10} $. 
By Lemmas~\ref{lem_type_fiber}(1) and \ref{lem_genus_k-2}, 
the genera of $F_{g \mathfrak{a}+ \mathfrak{b}}$ and $F_{g \mathfrak{a}+ 2\mathfrak{b}}$ are $g-2$ and $g$ respectively. 
Hence $\lambda_{(g,2)} =\min \Lambda_g(\tfrac{-3}{2}) $. 

One has 
$$\lambda_{(3,2)} \approx  1.50614>  \lambda_{(3,1)} = \lambda_{(4, 3)} \approx 1.40127,$$
and hence $\min \Lambda_3(\tfrac{-3}{2}) > \min \Lambda_3(\tfrac{-1}{2}) $. 
By Proposition~\ref{prop_monotone} together with the inequality 
$$\lambda_{(4,1)} \approx  1.28064>  \lambda_{(5,3)}  \approx 1.26123,$$
one obtains 
$\lambda_{(k,1)} > \lambda_{(k+1,3)}$ for any $k \ge 4$. 
The inequality 
$\min \Lambda_g(\tfrac{-3}{2}) > \min \Lambda_g(\tfrac{-1}{2}) $  holds for $g \equiv 3 \pmod{10}$ and $g >3$ 
since 
$$\min \Lambda_g(\tfrac{-3}{2}) = \lambda_{(g,2)} > \lambda_{(g,1)} > \lambda_{(g+1,3)} \ge \min \Lambda_g(\tfrac{-1}{2}).$$
One completes  the proof of (3).  
Similarly one can prove  (4). 
$\Box$

\begin{rem} 
\label{rem_our-example}
The pseudo-Anosov homeomorphism  whose dilatation equals $\min \Lambda_g(\tfrac{-3}{2})$ in the proof of 
Proposition~\ref{prop_bound_KT1}(1) (resp. (2)) is  non-orientable (resp. orientable), 
see Corollary~\ref{cor_kl_orientable}. 
\end{rem}

\noindent
{\it Proof of Theorem~\ref{thm_bound_KT1}.} 
See Proposition~\ref{prop_bound_KT1}(1),(2). 
$\Box$
\medskip

In the case $g \equiv 2,4 \pmod{10}$, we compute $\min \Lambda_g(\tfrac{-3}{2})$ under certain conditions of $g$. 

\begin{prop}
\label{prop_bound_KT2}
Let  $g \equiv 2,4 \pmod{10}$ and $g \ge 12$. 
Suppose that $g+2 \not \equiv 0$ $\pmod{4641 (= 3 \cdot 7 \cdot 13 \cdot 17)}$. 
\begin{description}
\item[(1)] 
$\lambda_{(g+2,3)} =\min \Lambda_g(\tfrac{-3}{2}) < \min \Lambda_g(\tfrac{-1}{2})$ if $\gcd(g+2,3)=1$. 
\item[(2)] 
$\lambda_{(g+2,7)} =\min \Lambda_g(\tfrac{-3}{2}) < \min \Lambda_g(\tfrac{-1}{2})$ if $3$ divides $g+2$ and $\gcd(g+2,7)=1$. 
\item[(3)] 
$\lambda_{(g+2,13)} =\min \Lambda_g(\tfrac{-3}{2}) < \min \Lambda_g(\tfrac{-1}{2})$ if $21(= 3 \cdot 7)$ divides $g+2$ and $\gcd(g+2,13)=1$. 
\item[(4)] 
$\lambda_{(g+2,17)} =\min \Lambda_g(\tfrac{-3}{2}) < \min \Lambda_g(\tfrac{-1}{2})$ if $273 (= 3 \cdot 7 \cdot 13)$ divides $g+2$ and $\gcd(g+2,17)=1$. 
\end{description}
\end{prop}

\noindent
The following will be used for proving  Proposition~\ref{prop_bound_KT2}. 
Its proof is similar to the one for Proposition~\ref{prop_bound_KT1}(3). 

\begin{lem}
\label{lem_step1}
\ 
\begin{description}
\item[(1)]  
Let $g \equiv 2 \pmod{10}$ and $g \ge 12$. 
\begin{description}
\item[(i)] 
Suppose that $g \equiv 2,22 \pmod{30}$. 
If $\lambda_{(g, \ell)} \in \Lambda_g(\tfrac{-3}{2})$, then $\ell \ge 3$. 
\item[(ii)]
Suppose that $g \equiv 12 \pmod{30}$. 
If $\lambda_{(g, \ell)} \in \Lambda_g(\tfrac{-3}{2})$, then $\ell \ge 5$. 
\end{description}
\item[(2)] 
Let $g \equiv 4 \pmod{10}$ and $g \ge 14$. 
Then $\lambda_{(g,1)} \in \Lambda_g(\tfrac{-3}{2})$. 
\end{description}
\end{lem}

\begin{lem}
\label{lem_step2}
Suppose that $g \equiv 2,4 \pmod{10}$ and $g \ge 12$. 
If $gcd(g+2,\ell)=1$, $\ell \equiv 2,3 \pmod{5}$ and $0 < \ell< g+2$, then 
$\lambda_{(g+2, \ell)} \in \Lambda_g(\tfrac{-3}{2})$. 
\end{lem}

\noindent
{\it Proof.} 
We use Lemma~\ref{lem_type_fiber}(1). 
We set $k= g+2$ ($k \equiv 1,4 \pmod{5}$). 
If $\ell \equiv 2,3 \pmod{5}$, then 
$\gcd(2k+ \ell, 5)= 5$ or $\gcd(5, k+ 2\ell) = 5$. 
Thus if $\ell$ satisfies that $\gcd(k, \ell) = \gcd(g+2, \ell)=1$ and $0 < \ell < g+2$, 
then one obtains the desired claim $\lambda_{(g+2, \ell)} \in \Lambda_g(\tfrac{-3}{2})$. 
$\Box$
\medskip

One can check the following inequalities. 
\begin{lem}
\label{lem_compare-min}
\begin{description}
%\item[(1)] 
%$ \lambda_{(5.3)} \approx 1.26123 < \lambda_{(4,1)} \approx 1.28064$. 
 \item[(1)] 
$\lambda_{(9,7)} \approx 1.16873< \lambda_{(8,1)} \approx 1.12876$. 
\item[(2)] 
$\lambda_{(73,13)} \approx 1.013457447 < \lambda_{(72,1)} \approx 1.013457858$. 
\item[(3)]
$\lambda_{(125,17)} \approx 1.007791640 < \lambda_{(124,1)} \approx 1.007791898$.  
\end{description}
\end{lem}

\noindent
{\it Proof of Proposition~\ref{prop_bound_KT2}.}
(1) 
By Lemma~\ref{lem_step2}, $\lambda_{(g+2,3)} \in \Lambda_g(\tfrac{-3}{2})$. 
We have shown that 
$\lambda_{(k+1,3)}< \lambda_{(k,1)}$ for any $k \ge 4$ in the proof of Proposition~\ref{prop_bound_KT1}(3),(4). 
Hence $\lambda_{(g+2,3)}< \lambda_{(g+1,1)}$ for any $g \ge 3$. 
By (\ref{equation_bound-one-half}), we have 
$\min \Lambda_g (\tfrac{-3}{2}) \le \lambda_{(g+2,3)} < \lambda_{(g+1,1)} \le  \min \Lambda_g (\tfrac{-1}{2})$. 
We can prove that $\lambda_{(g+2,3)}$ attains $\min \Lambda_g (\tfrac{-3}{2})$ 
by using the foregoing argument together with Lemma~\ref{lem_step1}(1). 

The claims (2),(3),(4) can be verified by using Lemmas~\ref{lem_step1}, \ref{lem_step2} and \ref{lem_compare-min}. 
$\Box$ 
\medskip

\noindent
{\it Proof of Theorem~\ref{thm_bound_KT2}.} 
See Proposition~\ref{prop_bound_KT2}. 
$\Box$

\begin{ques}
Is it true that 
$\delta_g \le \min \Lambda_g(\tfrac{-3}{2}) < \min \Lambda_g(\tfrac{-1}{2})$ for all $g \equiv 2,4 \pmod{10}$ and $g \ge 12$?
\end{ques}

\begin{rem}
\label{rem_AD}
Independently, 
Aaber and Dunfield identified the pair $(k(g), \ell(g))$ such that 
the pseudo-Anosov homeomorphism $\overline{\Phi}_{k(g) \frak{a}+ \ell(g) \frak{b}}: \varSigma_g \rightarrow \varSigma_g$ which 
attains $\min \Lambda_g(\tfrac{-3}{2})$ for large $g$. 
They proved that under a plausible assumption, 
the mapping class $\overline{\phi}_{k(g) \frak{a}+ \ell(g) \frak{b}} = [\overline{\Phi}_{k(g) \frak{a}+ \ell(g) \frak{b}}]$ has the least volume among 
pseudo-Anosov elements of $\mathrm{Mod}(\varSigma_g)$ for large $g$, see \cite{AD}. 
\end{rem}

We turn to  $ \min \Lambda_g^+(\tfrac{-3}{2})$. 
By Corollary~\ref{cor_kl_orientable}(1) and Lemma~\ref{lem_sym_fib}(1), 
one sees that if $g$ is even, then there exist no 
orientable pseudo-Anosov monodromies with a closed fiber of genus $g$ of fibrations on $N(\tfrac{-3}{2})$. 
Hence in this case $  \Lambda_g^+(\tfrac{-3}{2}) = \emptyset$. 
We compute $ \min \Lambda_g^+(\tfrac{-3}{2})$ for $g$ odd.

\begin{prop}
\label{prop_ori_bound_KT}
Let $g \ge 5$. 
\begin{description}
\item[(1)] 
$\lambda_{(g+2,2)}= \min \Lambda_g^+(\tfrac{-3}{2}) < \min \Lambda_g^+(\tfrac{-1}{2})$ if $g \equiv 7,9 \pmod{10}$. 
\item[(2)] 
$\lambda_{(g+2,4)}= \min \Lambda_g^+(\tfrac{-3}{2}) \le  \min \Lambda_g^+(\tfrac{-1}{2})$ if $g \equiv 1,5 \pmod{10}$. 
The equality holds if and only if $g=5$. 
\item[(3)] 
$\lambda_{(g,2)}= \min \Lambda_g^+(\tfrac{-3}{2}) > \min \Lambda_g^+(\tfrac{-1}{2})$ if $g \equiv 3 \pmod{10}$. 
\end{description}
\end{prop}

\noindent
{\it Proof.} 
We use Corollary~\ref{cor_kl_orientable} to see whether $\lambda_{(k,\ell)} \in \Lambda_g(\tfrac{-3}{2})$ is an element of $ \Lambda^+_g(\tfrac{-3}{2})$ or not. 
\medskip
\\
(1) 
We see that $\lambda_{(g+2,2)} \in \Lambda^+_g(\tfrac{-3}{2})$, see Remark~\ref{rem_our-example}. 
By Proposition~\ref{prop_bound_KT1}(2), we have 
$$\lambda_{(g+2,2)} = \min  \Lambda_g(\tfrac{-3}{2}) = \min \Lambda^+_g(\tfrac{-3}{2})< \min \Lambda_g(\tfrac{-1}{2}) \le  \min \Lambda^+_g(\tfrac{-1}{2}).$$
(2)
It can be shown that 
$\lambda_{(g+2,4)}  = \min \Lambda_g^+(\tfrac{-3}{2})$. 
Since $\lambda_{(7,4)}= \lambda_{(6,1)}$, the equality 
$ \min \Lambda_5^+(\tfrac{-3}{2}) =  \min \Lambda_5^+(\tfrac{-1}{2})$ holds. 
Suppose that $g \ne 5$. 
By Proposition~\ref{prop_monotone} together with 
$$\lambda_{(8,4)} \approx 1.14555 <  \lambda_{(7,1)} \approx 1.14879,$$
we obtain the inequality 
$\lambda_{(k,4)}< \lambda_{(k-1,1)}$ for any $k \ge 8$. 
Thus 
$ \min \Lambda_g^+(\tfrac{-3}{2}) <  \min \Lambda_g^+(\tfrac{-1}{2})$.  

One can prove (3) by using the similar argument together with Proposition~\ref{prop_bound_KT1}(3). 
$\Box$
\medskip

\noindent
{\it Proof of Theorem~\ref{thm_ori_bound_KT}.} 
See Proposition~\ref{prop_ori_bound_KT}(1),(2). 
$\Box$
\medskip

\noindent
{\it Proof of Proposition~\ref{prop_non-monotone}.} 
We have proved the inequality $(\lambda_{(g+2,2)}<) \lambda_{(g+2,4)}< \lambda_{(g+1,1)}$  for any $g  \ge 6$ 
in the proof of Proposition~\ref{prop_ori_bound_KT}(2). 
By Theorem~\ref{thm_ori_bound_KT} and by the assumption $ \delta_{g+1}^+ = \lambda_{(g+1,1)}$,  one has 
$$\delta_g^+ \le \max\{\lambda_{(g+2,2)}, \lambda_{(g+2,4)}\} \le  \lambda_{(g+2,4)}< \lambda_{(g+1,1)}= \delta_{g+1}^+.$$ 
This completes the proof. 
$\Box$

\begin{rem}\ 
\label{rem_-237}
\begin{description}
\item[(1)] 
The $(-2,3,7)$-pretzel knot complement is homeomorphic to $N(\tfrac{-3}{2}, \tfrac{-8}{3})$, see \cite[Table~A.4]{MP}. 
On the other hand, 
$\mathrm{slope}(7 \mathfrak{a}+ 4 \mathfrak{b})  = (\tfrac{-37}{18},  \tfrac{-3}{2}, \tfrac{-8}{3})$. 
The monodromy $\Phi_{7 \mathfrak{a}+ 4 \mathfrak{b}}: \varSigma_{5,17} \rightarrow  \varSigma_{5,17}$ 
of the fibration associated to  $7 \mathfrak{a}+ 4 \mathfrak{b}$ on $N$ is orientable (see Corollary~\ref{cor_kl_orientable}(1)) and its 
singularity data  is given by   $(16)$ (see Corollay~\ref{cor_sing_2cusps}(1)).  
Thus $\Phi_{7 \mathfrak{a}+ 4 \mathfrak{b}}: \varSigma_{5,17} \rightarrow \varSigma_{5,17}  $ extends to the pseudo-Anosov monodromy 
$\overline{\Phi}_{7 \mathfrak{a}+ 4 \mathfrak{b}}: \varSigma_{5,1} \rightarrow \varSigma_{5,1}$ of the fibration 
on $N(\tfrac{-3}{2}, \tfrac{-8}{3})$ (with the dilatation $\lambda_{(7,4)}$) 
by capping all the boundary components on $T_{\beta} \cup T_{\gamma}$. 
\item[(2)] 
$\overline{\Phi}_{7 \mathfrak{a}+ 4 \mathfrak{b}}: \varSigma_{5,1} \rightarrow \varSigma_{5,1}$ 
extends to  the  monodromy $: \varSigma_{5} \rightarrow \varSigma_{5} $ 
of the fibration on $N(\tfrac{-37}{18},  \tfrac{-3}{2}, \tfrac{-8}{3})$ 
with dilatation  $\delta_5^+ =\lambda_{(7,4)}$. 
Since this extended monodromy is orientable, we have  $\delta_5^+ \in \Lambda_5^+(\tfrac{-3}{2})$. 
\end{description}
\end{rem}

\noindent
By Remark~\ref{rem_-237}(2), we have: 

\begin{lem}
\label{lem_list-KT2}
$N(\tfrac{-3}{2}) \in \mathcal{U}^+$. 
\end{lem}

\subsection{Fibers of genera $8$ and $13$}

By using the foregoing discussion one can prove the following which implies Proposition~\ref{prop_small_genus}. 

\begin{prop}\ 
\begin{description}
\item[(1)] 
$N(\tfrac{-4}{3}, \tfrac{-25}{17}, -5)$ is a $\varSigma_8$-bundle over the circle 
with dilatation $\lambda_{(18,17,7)} \thickapprox 1.10403$ and with the singularity data $(\underbrace{1, \cdots, 1}_6, 15, \underbrace{1, \cdots, 1}_7)$. 
\item[(2)] 
$N(\tfrac{-29}{27}, \tfrac{-5}{3}, -6)$ is a $\varSigma_{13}$-bundle over the circle 
with dilatation $\lambda_{(27,21,8)} \thickapprox 1.07169$ and with  the singularity data $(25, \underbrace{1, \cdots, 1}_7, \underbrace{2, \cdots, 2}_8)$. 
\end{description}
\end{prop}

\noindent
Department of Mathematical and Computing Sciences, Tokyo Institute of Technology \\
  Ohokayama, Meguro  Tokyo 152-8552 Japan \\
 E-mail address: kin@is.titech.ac.jp
 \medskip
 
 \noindent
 Department of Mathematical and Computing Sciences,  Tokyo Institute of Technology \\
  Ohokayama, Meguro  Tokyo 152-8552 Japan \\
 E-mail address: takasawa@is.titech.ac.jp

\end{document}